%% file: main.tex
\documentclass[12pt]{amsart}
\pdfoutput=1

\usepackage[utf8]{inputenc}
\usepackage[margin=1in]{geometry}
\usepackage{graphicx, url} 
\usepackage{amsmath,amssymb,amsthm}
\usepackage{color, tikz}
\usepackage{tabularx}
\usepackage[hidelinks]{hyperref}

\usepackage{bbm}

\usepackage{rotating}
\usepackage{tikz,graphicx}
\tikzstyle{vertex}=[circle, draw, inner sep=0pt, minimum size=4pt]
\newcommand{\vertex}{\node[vertex]}

\theoremstyle{plain}
\newtheorem{theorem}{Theorem}[section]
\newtheorem{prop}[theorem]{Proposition}

\newtheorem{lemma}[theorem]{Lemma}
\newtheorem{corollary}[theorem]{Corollary}

\theoremstyle{definition}
\newtheorem{remark}[theorem]{Remark}
\newtheorem{definition}[theorem]{Definition}

\newtheorem{example}[theorem]{Example}

\newcommand{\indeg}{\mathrm{indeg}}
\newcommand{\outdeg}{\mathrm{outdeg}}
\newcommand{\rsimplex}{\Delta_\mathcal{R}}
\newcommand{\eq}{\mathrm{eq}}
\newcommand{\dkkeq}{T_{\eq}}
\newcommand{\flow}{\mathcal{F}_1}
\newcommand{\op}[1]{\mathcal{O}(#1)}
\newcommand{\Out}{\mathrm{Out}}
\newcommand{\In}{\mathrm{In}}
\newcommand{\cone}[1]{\mathrm{cone}(#1)}
\newcommand{\inlev}{\mathrm{inlevel}}
\newcommand{\outlev}{\mathrm{outlevel}}

\newcommand{\conv}[1]{\mathrm{conv}\{#1\}}
\newcommand{\R}{\mathbb{R}}

\newcommand{\Z}{\mathbb{Z}}

\newcommand{\Or}{\mathcal{O}}

\definecolor{green}{RGB}{34, 139, 34}

\newcommand\commentout[1]{}

\DeclareFontFamily{U}{mathx}{}
\DeclareFontShape{U}{mathx}{m}{n}{<-> mathx10}{}
\DeclareSymbolFont{mathx}{U}{mathx}{m}{n}
\DeclareMathAccent{\widecheck}{0}{mathx}{"71}

\title{Equatorial Flow Triangulations of Gorenstein Flow Polytopes}

\author{Benjamin Braun}
\address{Department of Mathematics\\
         University of Kentucky\\
\url{https://sites.google.com/view/braunmath/}}
\email{benjamin.braun@uky.edu}

\author{Alvaro Cornejo}
\address{Department of Mathematics\\
         University of Kentucky\\
\url{https://sites.google.com/view/alvaro-cornejo/home}}
\email{alvaro.cornejo@uky.edu}

\date{14 February 2025}

\thanks{Both authors were partially supported by the US National Science Foundation award DMS-1953785.
The authors thank Jonah Berggren for sharing his preliminary results regarding projection maps of Gorenstein flow polytopes for amply framed DAGs, which inspired the contents of Section~\ref{sec:projection}.}

\begin{document}

\begin{abstract}
Generalizing work of Athanasiadis for the Birkhoff polytope and Reiner and Welker for order polytopes, in 2007 Bruns and R\"omer proved that any Gorenstein lattice polytope with a regular unimodular triangulation admits a regular unimodular triangulation that is the join of a special simplex with a triangulated sphere.
These are sometimes referred to as equatorial triangulations.
We apply these techniques to give purely combinatorial descriptions of previously-unstudied triangulations of Gorensten flow polytopes.
Further, we prove that the resulting equatorial flow polytope triangulations are usually distinct from the family of triangulations obtained by Danilov, Karzanov, and Koshevoy via framings.
We find the facet description of the reflexive polytope obtained by projecting a Gorenstein flow polytope along a special simplex.
Finally, we show that when a partially ordered set is strongly planar, equatorial triangulations of a related flow polytope can be used to produce new unimodular triangulations of the corresponding order polytope.
\end{abstract}

\maketitle

\section{Introduction}

Given a finite directed acyclic graph (DAG) $G$, the \emph{flow polytope} $\flow(G)$ is the lattice polytope of non-negative flows on $G$ of strength one.
Flow polytopes have been the subject of intense recent study in geometric and algebraic combinatorics~\cite{baldoni-vergne,caracolvolume,combflowmodel,braunmcelroyfullvolumes,DKK,generalpitmanstanley,cyclicorderflow,integerpointsflow,lidskii,gtflow,flowdiagonalharmonics,meszaros-morales,meszaros-morales-striker,morrisidentityflow,rietsch2024rootpolytopesflowpolytopes,vonbell2024triangulations,framedtriangulations}.
Flow polytopes are connected to many areas of mathematics; one example relevant to this work is that when the DAG $G$ is strongly planar, it is known that the flow polytope of $G$ is integrally equivalent to the order polytope of the poset with Hasse diagram given by the strongly planar dual of $G$~\cite{meszaros-morales-striker}.

A major tool used in the study of flow polytopes is the \emph{DKK triangulation} induced by a framing on the edges of $G$, a regular unimodular triangulation introduced by Danilov, Karzanov, and Koshevoy~\cite{DKK}.
M\'esz\'aros, Morales, and Striker proved that when a strongly planar DAG $G$ is equipped with the planar framing, the induced DKK triangulation is integrally equivalent to the canonical triangulation of the order polytope for the poset dual to $G$.
DKK triangulations have been used extensively in the literature to study the structure of flow polytopes, dual graphs of their triangulations, and their volumes.

In this paper, we describe \emph{equatorial flow triangulations}, a new family of regular unimodular triangulations for Gorenstein flow polytopes that arise as the join of a unimodular simplex and a simplicial sphere.
The study of triangulations with this structure originated in work of Athanasiadis~\cite{athanasiadis} for the Birkhoff polytope and Reiner and Welker~\cite{reinerwelker} for order polytopes.
In 2007, Bruns and R\"omer~\cite{brunsromer} proved that any Gorenstein lattice polytope with a regular unimodular triangulation admits a regular unimodular triangulation that is the join of a special simplex with a triangulated sphere.
We apply this Bruns-R\"omer machinery to DKK triangulations arising from route decompositions of DAGs, producing equatorial flow triangulations of Gorenstein flow polytopes.

Our main contributions in this paper are the following.
\begin{enumerate}
\item In Theorem~\ref{thm:routedecomp}, we use route decompositions to give a new combinatorial characterization of DAGs yielding Gorenstein flow polytopes.
\item In Theorem~\ref{thm:conditionsflowequatorial}, we give a purely combinatorial description of equatorial flow triangulations of a Gorenstein flow polytope.
\item In Theorem~\ref{thm:notfullnotdkk}, we prove that if $G$ has an inner vertex with indegree at least $3$, then the equatorial flow triangulation is not a DKK triangulation.
\item As a corollary of the work of Bruns and R\"omer~\cite{brunsromer}, a Gorenstein $\flow(G)$ projects onto a reflexive polytope. 
In Section~\ref{sec:projection}, we provide a combinatorial halfspace description of this reflexive polytope for an equatorial flow triangulation.
\item Finally, in Section~\ref{sec:order}, we prove that when $G$ is strongly planar, a particular choice of route decomposition yields an equatorial flow triangulation that is integrally equivalent to the Reiner-Welker equatorial triangulation of the associated order polytope.
We also use equatorial flow triangulations to produce new triangulations of Gorenstein order polytopes for strongly planar posets.
\end{enumerate}

We begin with Section~\ref{sec:background}, which contains background on polytopes, triangulations, and Ehrhart theory.

\section{Background}\label{sec:background}

\subsection{Lattice and Flow polytopes}

Let $G$ be a finite directed acyclic graph (DAG) with linearly ordered vertex set $\{s<1<2<\cdots <n<t\}$ such that if $(i,j)$ is a directed edge in $G$, then $i<j$.
We assume throughout this work that $s$ is the unique source of $G$ and $t$ is the unique sink of $G$, and all other vertices of $G$ we call \emph{inner vertices}.
We denote by $\mathrm{in}_G(v)$ the set of incoming edges at $v$, and by $\mathrm{out}_G(v)$ the set of outgoing edges at $v$.
A \emph{route} is a directed path in $G$ from $s$ to $t$.

\begin{definition}\label{def:flowpolytope}
    The \emph{flow polytope} of $G$ is
    \[
    \flow(G):=\left\{x\in \R_{\geq 0}^{E(G)}: \sum_{e\in \mathrm{in}_G(v)}x_e=\sum_{e\in \mathrm{out}_G(v)}x_e \text{ for every inner vertex }v, \sum_{e\in \mathrm{out}_G(s)}x_e=1 \right\} \, .
    \]
    Equivalently, $\flow(G)$ is given by the convex hull of indicator vectors of routes in $G$.
    The \emph{cone of non-negative flows} of $G$ is 
     \[     \mathcal{F}(G):=\left\{x\in \R_{\geq 0}^{E(G)}: \sum_{e\in \mathrm{in}_G(v)}x_e=\sum_{e\in \mathrm{out}_G(v)}x_e \text{ for every inner vertex }v \right\} \, .
    \]
\end{definition}

The dimension of $\flow(G)$ is $|E|-n-1$, where $G$ has $n$ inner vertices and $E$ is the set of edges of $G$. 
Note that the facets of $\flow(G)$ and $\mathcal{F}(G)$ are given by $x_e = 0$ for an edge $e \in E(G)$, but not all of these equations are facet defining. 
Towards identifying the facet defining hyperplanes, we define an edge to be \emph{idle} if it is the only incoming or outgoing edge edge of an inner vertex. 
It is known~\cite{DKK,vonbell2024triangulations} that contracting an idle edge $e \in E(G)$ does not change the lattice-polyhedral structure of $\flow(G)$ and $\mathcal{F}(G)$.
Further, in a DAG with no idle edges, every edge yields a facet hyperplane $x_e=0$. 
Hence, we will assume throughout this paper that we have no idle edges, and this can be done without losing information about $\flow(G)$ and $\mathcal{F}(G)$.

A \emph{framing} of $G$ is a choice for each inner vertex $v$ of linear orderings on $\mathrm{in}_G(v)$ and $\mathrm{out}_G(v)$.
Danilov, Karzanov, and Koshevoy~\cite{DKK} proved that every framing $F$ induces a regular unimodular triangulation of $\flow(G)$, called the \emph{DKK triangulation} induced by $F$, which we explain next.
Assume that we have chosen a framing $F$ of $G$, and note that all of the following definitions are dependent on $F$.
If $P$ is a route in $G$ containing the vertex $v$, then we write $Pv$ for the path from the source $s$ to $v$ following $P$, and similarly $vP$ for the path from $v$ to the sink $t$ along $P$.
We write $vPw$ for the segment of $P$ from $v$ to $w$.
We write $\Out(v)$ for the set of partial routes in $G$ starting at $v$ and ending at $t$, and similarly $\In(v)$ is partial routes from $s$ to $v$.
If $e$ is less than $f$ in the linear order for $F$ on $\mathrm{in}(v)$, we write $e\prec_{\mathrm{in}(v)} f$, and similarly for $\mathrm{out}(v)$.

Let $P$ and $Q$ be paths in $\Out(v)$ that agree on the paths $P'\subset P$ and $Q' \subset Q$ that begin at $v$ and end at $w$.
Suppose further that the vertices following $w$ on $P$ and $Q$ are distinct, and call them $w_P$ and $w_Q$.
We define an inequality between these paths by setting $P\prec_{\Out(v)} Q$ if $(w,w_P)\prec_{F,\mathrm{out}(w)} (w,w_Q)$, and similarly for $P\prec_{\In(v)} Q$.
Assume $P$ and $Q$ are routes that intersect at a common inner vertex $v$ of $G$.
We say $P$ and $Q$ are \emph{in conflict}, also called \emph{conflicting}, if $Pv\prec_{\In(v)} Qv$ and $vQ \prec_{\Out(v)} vP$.
If $P$ and $Q$ are not conflicting at $v$, then they are \emph{coherent} at $v$.
$P$ and $Q$ are called \emph{coherent} if they are coherent at every inner vertex $v$ that is contained in both $P$ and $Q$.
Finally, a \emph{clique} is a set of pairwise-coherent routes in $G$.
When a route $R$ in $G$ is coherent with every other route in $G$, $R$ is called \emph{exceptional}.

\begin{definition}\label{def:dkktriangulation}
Given a DAG $G$ with framing $F$, the \emph{DKK triangulation} of $\flow(G)$ induced by $F$ has one facet for each maximal clique, where the facet for a clique is the simplex given by the convex hull of the indicator vectors for the routes in the clique.
\end{definition}

Recall that a lattice simplex $\conv{v_0,v_1,\ldots,v_d}\subset \R^d$ is \emph{unimodular} if $\{v_1-v_0,\ldots,v_d-v_0\}$ is a lattice basis for $\Z^d$.
Further, a triangulation of a lattice polytope is \emph{unimodular} if every simplex in the triangulation is unimodular.

\begin{theorem}[Danilov, Karzanov, and Koshevoy~\cite{DKK}]\label{thm:dkk}
    For any finite DAG $G$ and any framing $F$, the DKK triangulation of $\flow(G)$ induced by $F$ is a regular unimodular triangulation.
\end{theorem}

\begin{figure}
    \centering
    \input{clique.tikz}
    \caption{An example of a framed DAG where the edges are linearly ordered by $1 < 2$ and an example of coherent routes forming a maximal clique. The bottom two routes are the exceptional routes for this framed DAG.}
    \label{fig:clique}
\end{figure}
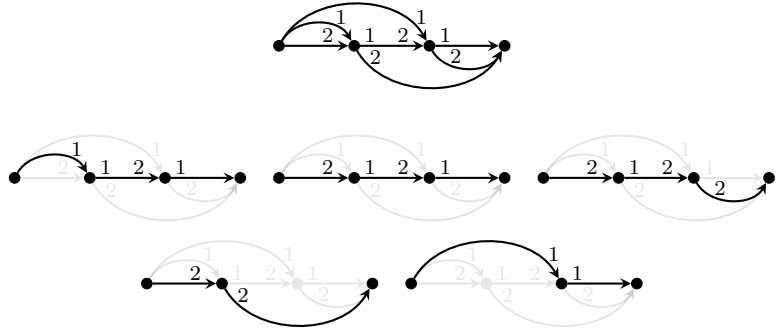

\begin{example}\label{ex:clique}
    Figure~\ref{fig:clique} is an example of a framed DAG and coherent routes.
\end{example}

\subsection{Face Enumeration and Ehrhart Theory}

We next recall $h$-polynomials, Ehrhart series, and connections among these objects.
Two invariants of interest for flow polytopes and their triangulations are $h$-vectors and Ehrhart $h^*$-vectors, defined as follows.

\begin{definition}
\label{def:hpoly}
    Given a finite simplicial complex $T$, let $f_k$ denote the number of $k$-dimensional faces of $T$, where we define $f_{-1}:=1$.
The \emph{$h$-polynomial} of $T$ is defined to be
\[
h(T;z):=\sum_{k=-1}^{d-1}f_kz^{k+1}(1-z)^{d-1-k} \, .
\]
The coefficient vector of $h(T;z)$ is the \emph{$h$-vector} of $T$.
\end{definition}

In the following definition, we rely on a result due to Ehrhart~\cite{Ehrhart} that states the Ehrhart series of a lattice polytope is a rational function of the form given in the definition.

\begin{definition}
\label{def:ehrhart}
    Given polytope $P$ of dimension $d$ with vertices in $\Z^n$, the \emph{Ehrhart series} of $P$ is
    \[
    \mathrm{Ehr}(P;z):=1+\sum_{t\geq 1}|tP\cap \Z^d|z^t = \frac{h_0^*+h_1^*z+\cdots + h_d^*z^d}{(1-z)^{d+1}}\, .
    \]
    We call the numerator of this rational function the \emph{$h^*$-polynomial} of $P$ and its coefficient vector the \emph{$h^*$-vector} of $P$.
\end{definition}

\begin{example}
\label{ex:ehrhartforflow}
The Ehrhart series for the flow polytope of the DAG $G$ given in Figure~\ref{fig:clique} is 
\[
    \mathrm{Ehr}(\flow(G);z) = \frac{1+3z+z^2}{(1-z)^5}
\]
\end{example}

For a unimodular triangulation, the $h$-vector and $h^*$-vector are closely related; see Beck and Robins~\cite{ccd} for a textbook proof of the following result.

\begin{theorem}\label{thm:hhstar}
Given a unimodular triangulation $T$ of a lattice polytope $P$, the $h$-vector of $T$ is equal to the $h^*$-vector of $P$.
\end{theorem}

We complete this subsection by recalling the definition of the codegree of a lattice polytope.

\begin{definition}
    The \emph{degree} of a $d$-polytope $P$ is 
    \[
        \deg(P) := \max\{m \in \mathbb{Z}_{\geq 0} : h_m^* \neq 0\}  \, .
    \]

    The \emph{codegree} of $P$ is
    \[
        \mathrm{codeg}(P) := d+1 - \deg(P) \, .
    \]
\end{definition}

The following theorem is a consequence of Ehrhart-Macdonald reciprocity~\cite[Theorem~4.5]{ccd}.

\begin{theorem}\label{thm:codegreeinterior}
    The codegree of a $d$-polytope $P$ is the smallest integer $k$ such that $kP$ contains an interior lattice point. 
\end{theorem}

\begin{example}
    The degree and codegree for the flow polytope of the DAG given in Figure~\ref{fig:clique} is $\deg(\flow(G)) = 2$ and $\mathrm{codeg}(\flow(G)) = 4 + 1 - 2 = 3$. Hence, the dilate $3 \flow(G)$ contains an interior lattice point by Theorem~\ref{thm:codegreeinterior}.
\end{example}

\subsection{Gorenstein Polytopes and Equatorial Triangulations}

We next recall Gorenstein polytopes and their properties.

\begin{definition}\label{def:gorcone}
    A rational pointed cone $C\subseteq \mathbb{R}^{d+1}$ is said to be \emph{Gorenstein} if there exists an integer point $c \in \mathbb{Z}^{d+1}$ such that
    \[
        c + (C \cap \mathbb{Z}^{d+1}) = C^\circ \cap \mathbb{Z}^{d+1}
    \]
    where $C^\circ$ is the (relative) interior of $C$. 
    We call $c$ the \emph{Gorenstein point} of $C$. 
\end{definition}

We can extend the definition of Gorenstein to lattice polytopes as follows.

\begin{definition}\label{def:gorpolytope}
    A lattice polytope $P \subseteq \mathbb{R}^d$ is \emph{Gorenstein} if \[
        \cone{P} := \{(x,t) \in \mathbb{R}^{d+1} : t \geq 0, x \in tP \}
    \]
    is a Gorenstein cone.
    The polytope $P$ is \emph{reflexive} if the Gorenstein point in $\cone{P}$ has final coordinate $1$, i.e., is at height $1$ in the cone over $P$.
\end{definition}

Note that if $P$ is $d$-dimensional, then $P$ is reflexive if and only if $P$ contains a unique interior lattice point $v$ and there exists an integer matrix $A$ such that 
\[
P-v = \{x\in \R^d: Ax\leq \mathbbm{1}\} \, .
\]
The following two propositions connect arbitrary Gorenstein polytopes to reflexive polytopes; for more details see the textbook treatment by Bruns and Herzog~\cite{BrunsHerzogCMR}.

\begin{prop}\label{prop:gorensteinreflexivecondition}
    A lattice polytope $P$ is Gorenstein if and only if there exists some $r \in \mathbb{N}$ such that $rP$ is reflexive.    
    In this case, we say $P$ is Gorenstein of index $r$ and, in addition, $r$ is the codegree of $P$.
\end{prop}

\begin{prop}\label{prop:gorensteinpointindex}
    If $P \subseteq \mathbb{R}^d$ is a Gorenstein polytope with Gorenstein point $c \in \mathbb{R}^{d+1}$ for $\cone{P}$, then $c = (v,r)$ where $r \in \mathbb{Z}_{\geq 1}$ is the index of $P$ and $v$ is the unique interior lattice point of $rP$.
\end{prop}

A key idea in the theory of equatorial triangulations is that of a special simplex, defined as follows.

\begin{definition}\label{def:specialsimplex}
    Given a $d$-dimensional lattice polytope $P$, a simplex $\Delta = \conv{v_1,\ldots,v_k}$ with $v_i \in P \cap \Z^d$ is \emph{special} if $\Delta \cap F$ is a facet of $\Delta$ for all facets $F$ of $P$.
\end{definition}

For Gorenstein lattice polytopes, the following proposition provides a mechanism for finding special simplices.

\begin{prop}[Bruns and R\"omer~\cite{brunsromer}]\label{prop:specialsimplex}
    If $P$ is a Gorenstein polytope and there exist lattice points $\{v_1,\ldots,v_k\} \subseteq P \cap \mathbb{Z}^d$ where $\{(v_1,1),\ldots,(v_k,1)\}$ sum to the Gorenstein point of $\cone{P}$, then $\conv{v_1,\ldots,v_k}$ forms a special simplex and is unimodular.
\end{prop}

The notion of special simplices were first introduced by Athanasiadis \cite{athanasiadis} and later explored by Bruns and R\"omer \cite{brunsromer} to show that Gorenstein polytopes with certain properties have unimodal $h^*$ coefficients. 
In particular Bruns and R\"omer proved that Gorenstein polytopes with a regular, unimodular triangulation have such unimodality. 
In this setting, they abstractly define a triangulation given as a join of a special simplex and the boundary of a simplicial polytope, as stated in the results below.

\begin{definition}\label{def:equatorialcomplexbackground}
    Given a polytope $P$ with a special simplex $\Delta=\conv{v_1,\ldots,v_k}$, define $\Gamma(P,\Delta)$ to be the polyhedral subcomplex of $\partial P$ generated by faces of the form $\bigcap_{i=1}^k F_i$ where $F_i$ is a facet of $P$ such that $v_i \notin F_i$. We call $\Gamma(P,\Delta)$ the \emph{equatorial complex} of $P$ with respect to $\Delta$, or just equatorial complex if $P$ and $\Delta$ are clear from context.
\end{definition}

\begin{theorem}[Bruns and R\"omer~\cite{brunsromer}]\label{thm:gortriangulations}
    If $\Delta \subseteq P$ is a special simplex of a Gorenstein polytope and $S$ is a triangulation of $\Gamma(P,\Delta)$, then $P$ is triangulated by the join $S \ast \Delta$, i.e., by the triangulation generated by $\conv{F \cup \Delta}$ where $F$ is a face of the triangulation $S$. Moreover, $S \ast \Delta$ is unimodular if $S$ is unimodular and $S$ is regular if it is the restriction to $\Gamma(P,\Delta)$ of a regular triangulation.
\end{theorem}

\begin{theorem}[Bruns and R\"omer~\cite{brunsromer}]\label{thm:gorrestrictionsimplicial}
    If $\Delta \subseteq P$ is a special simplex of a Gorenstein polytope $P$ and $S$ a regular unimodular triangulation of $P$, then the restriction $S \vert_{\Gamma(P,\Delta)}$ is the boundary complex of a simplicial polytope. 
\end{theorem}

Bruns and R\"omer also observed the following corollary.

\begin{corollary}\label{cor:gorensteinhhstar}
    If $\Delta \subseteq P$ is a special simplex of a Gorenstein polytope $P$ and $S$ a unimodular triangulation of $\Gamma(P,\Delta)$, then the $h^*$-polynomial of $P$ and $h$-polynomial of $S$ agree. 
    In other words,
    \[
        h_0^* + h_1^* z + \cdots + h_d^* z^d = h(S;z) .
    \]
\end{corollary}


\section{Gorenstein Flow Polytopes and Route Decompositions}

In this section, we provide a new characterization of DAGs with Gorenstein flow polytopes. Recall, that we assume $G$ has no idle edges.

\begin{definition}\label{def:degreeequality}
    A DAG $G$ is said to satisfy \emph{degree equality} if $\indeg_G(v) = \outdeg_G(v)$ for all inner vertices $v$.
\end{definition}

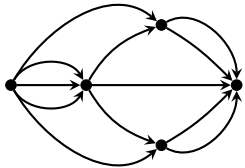
\begin{figure}[h]
    \centering
    \input{degreeequalityDAG.tikz}
    \caption{Example of a DAG with degree equality.}
    \label{fig:degreeequality}
\end{figure}

The following result gives a characterization of DAGs having a Gorenstein flow polytope.

\begin{theorem}[von Bell et al. ~\cite{vonbell2024triangulations}]\label{thm:degreeequality}
    The flow polytope $\flow(G)$ is Gorenstein if and only if $G$ satisfies degree equality.
\end{theorem}

The following corollary shows that every flow polytope can be found as a face of a Gorenstein flow polytope.

\begin{corollary}\label{cor:gorensteinface}
    Every flow polytope is a face of a Gorenstein flow polytope. 
\end{corollary}

\begin{proof}
Let $G$ be a DAG with flow polytope $\flow(G)$.
For every inner vertex $v$ of $G$, add edges either from the source $s$ to $v$ or from $v$ to the sink $t$ in order to create a DAG $G'$ containing $G$ that satisfies degree equality.
Then, $\flow(G)$ is obtained as a face of $\flow(G')$ by setting all the flow values on the newly added edges equal to zero.
\end{proof}

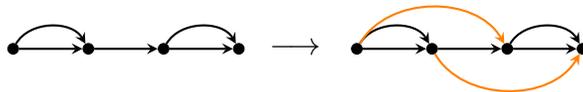
\begin{figure}[h]
    \centering
    \input{gorensteinface.tikz}
    \caption{Example of adding new edges (in orange) to a DAG to satisfy degree equality.}
    \label{ex:gorensteinface}
\end{figure}

Our goal in this section is to provide an alternative characterization of these DAGs using route decompositions, defined as follows.

\begin{definition}\label{def:routedecomp}
    A \emph{route decomposition} $\mathcal{R}$ of a DAG $G$ is a set of routes in $G$ such that the edge set of $G$ is the disjoint union of the edges of routes in $\mathcal{R}$.
\end{definition}

\begin{figure}
    \centering
    \input{routedecomp.tikz}
    \caption{Example of a route decomposition for the DAG found in figure \ref{fig:degreeequality}.}
    \label{fig:routedecomp}
\end{figure}
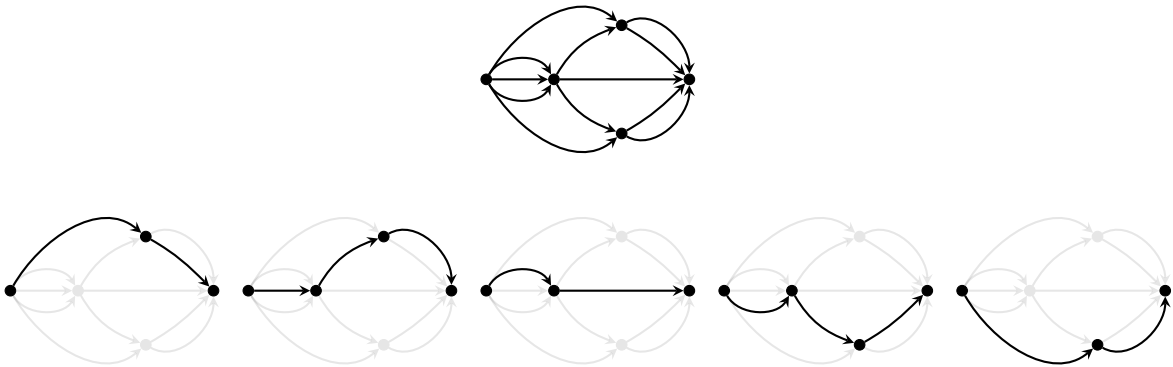

\begin{theorem}\label{thm:routedecomp}
    A DAG $G$ satisfies degree equality if and only if $G$ admits a route decomposition.
\end{theorem}

\begin{proof}
    The forward implication is proved in Proposition~\ref{prop:degreegivesdecomp} below.
    For the reverse implication, let $v \in V(G)$ be an inner vertex.
    As $v$ is an inner vertex, it has incoming and outgoing edges and thus must be incident to at least one route in our decomposition; let $\mathcal{S} \subseteq \mathcal{R}$ be the subset of our route decomposition containing routes incident to $v$. 
    For each $e \in \mathrm{in}_G(v)$, there exists a unique route $R_e \in \mathcal{S}$ since each edge is in some route of our decomposition and we have uniqueness as routes in $\mathcal{R}$ are edge disjoint. 
    Similarly, given $R \in \mathcal{S}$, there exists a unique $e \in \mathrm{in}_G(v)$ contained in $R$ since $G$ is acyclic. 
    Thus, we have a bijection between edges in $\mathrm{in}_G(v)$ and routes in $\mathcal{S}$ incident to $v$.
    We can make the same argument for $e' \in \mathrm{out}_G(v)$ and therefore
    \[
        \outdeg_G(v) = |\mathcal{S}| = \indeg_G(v) \, ,
    \]
    which proves that degree equality.
\end{proof}

The remainder of this section develops the content required to prove Proposition~\ref{prop:degreegivesdecomp}, which establishes the forward direction of Theorem~\ref{thm:routedecomp}.
For a DAG $G$ and a route $R$ from the source to the sink, we denote by $G-R$ the DAG with edge set $E(G)-E(R)$ and vertex set obtained by deleting from $V(G)$ any internal vertices that are incident only to edges in $R$.

\begin{figure}[h]
    \centering
    \input{deleteroute.tikz}
    \caption{An example of deleting routes}
    \label{fig:deleteroute}
\end{figure}
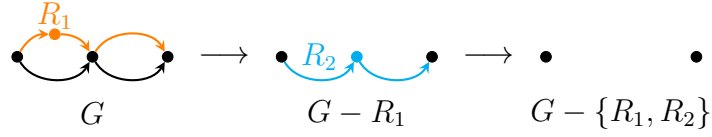

\begin{prop}\label{prop:stillacyclic}
    For a DAG $G$ and a route $R$, the DAG $G-R$ is acyclic.
\end{prop}

\begin{proof}
    If $G-R$ contains a cycle, then these edges are also in $G$ and thus $G$ contains a cycle, yielding a contradiction. 
\end{proof}

\begin{prop}\label{prop:deleterouteequality}
    If $G$ is a DAG satisfying degree equality, then $G-R$ also satisfies degree equality.
\end{prop}

\begin{proof}
    Consider $v \in V(G-R)$ an inner vertex of $G-R$.
    Since $v$ is not a source nor a sink in $G-R$, we have that $v$ must also be an inner vertex of $G$.
   Thus, we have degree equality for $v$ in $G$. 
    Now consider the following two cases. 
    If $v \notin V(R)$, then 
    \[
        \indeg_{G-R}(v) = \indeg_G(v) = \outdeg_G(v) = \outdeg_{G-R}(v) \, ,
    \]
    hence we have degree equality of $G-R$ at $v$. 
    Alternatively, if $v \in V(R)$, then the route $R$ must use exactly one incoming edge and exactly one outgoing edge of $v$ in $G$, since $G$ is acyclic. 
    Thus,
    \[
        \indeg_{G-R}(v) = \indeg_G(v) - 1 = \outdeg_G(v) - 1 = \outdeg_{G-R}(v) \, ,
    \]
    showing we have degree equality of $G-R$ at $v$. 
    This proves that $G-R$ satisfies degree equality. 
\end{proof}

\begin{lemma}\label{lem:innervertexroute}
    Suppose $G$ is a DAG that satisfies degree equality. 
    Then for any route $R$ in $G$, and for any inner vertex $v$ of $G-R$, there is a route in $G - R$ from $s$ to $t$ that is incident to $v$.
\end{lemma}
\begin{proof}
    Let $v$ be an inner vertex of $G-R$. 
    Thus, either $\indeg_{G-R}(v) \geq 1$ or $\outdeg_{G-R}(v) \geq 1$.
    Degree equality of $G-R$ follows from Proposition~\ref{prop:deleterouteequality}, and thus we have $\indeg_{G-R}(v) = \outdeg_{G-R}(v) \geq 1$. So there exists vertices $v_1,w_1 \in V(G-R)$ such that there are directed edges $(v_1,v)$ and $(v,w_1)$ in $G-R$. Consider $v_1$. If $v_1 = s$, then we have a directed path from $s$ to $v$. Otherwise, $v_1$ is an internal vertex of $G$ and degree equality of $G-R$ implies there exists a $v_2 \in V(G-R)$ such that we have a directed edge from $(v_2,v_1)$ in $G-R$. 
    We continue this process, constructing a sequence of $v_i \in V(G-R)$.
    We claim each $v_i$ is distinct. 
    If not, then we would have a directed cycle in $G-R$, which is a contradiction. 
    Since $G$ is finite, this process must end at $s$. 
    The same argument holds to establish a directed path from $w_j$ to $t$.
    Concatenating these two directed paths yields a route from $s$ to $t$ that is incident to $v$.
\end{proof}

\begin{corollary}\label{cor:internalconnected}
    Suppose $G$ is a DAG that satisfies degree equality. 
    For any route $R$ in $G$, if $G-R$ has internal vertices, the underlying undirected graph is connected.
\end{corollary}

\begin{lemma}\label{lem:equalsroute}
    For $G$ a DAG that satisfies degree equality and has one source and one sink, $G-R$ is disconnected if and only if $G=R$.
\end{lemma}

\begin{proof}
    If $E(G-R) = \emptyset$, then by definition $G-R$ is disconnected as it consists of only the source and sink with no edges, establishing the reverse implication. 
    For the forward implication, suppose $G-R$ is disconnected. By Corollary~\ref{cor:internalconnected}, we see $G-R$ has no internal vertices. 
    Thus, any edges in $G$ were in $R$, along with any internal vertices in $G$.
\end{proof}

\begin{lemma}\label{lem:samesink}
    For $G$ a DAG that satisfies degree equality, if $G-R$ is connected then it has the same unique source and sink.
\end{lemma}

\begin{proof}
    Suppose we have a second source $v \in V(G-R)$, i.e., we have $v \neq s$ and $\indeg_{G-R}(v) = 0$. 
    We claim $v \neq t$. 
    Suppose for contradiction that $v=t$, so we have $\indeg_{G-R}(t) = 0$ and we have
    \[
        \outdeg_{G-R}(t) = \outdeg_G(t) = 0 \, .
    \]
    It follows that $G-R$ is disconnected which is a contradiction.
    Hence $v$ is an inner vertex of $G$ and satisfies $\outdeg_{G-R}(v) = \indeg_{G-R}(v) = 0$. 
    This implies that $G-R$ is disconnected, but again by our assumption that $G - R$ is connected, this is a contradiction. 
    Hence no such $v$ can exist. 
    A similar argument shows that there is a unique sink.
\end{proof}

\begin{prop}\label{prop:degreegivesdecomp}
    Every DAG $G$ with one source and one sink that satisfies degree equality admits a route decomposition.
\end{prop}

\begin{proof}
    Let $R_1$ be a route from the source to the sink in $G$.
    We have $G-R_1$ is acyclic by Proposition~\ref{prop:stillacyclic} and satisfies degree equality by Proposition~\ref{prop:deleterouteequality}. 
    If $G-R_1 = \emptyset$ then we have a route decomposition. 
    Otherwise, we have the same source $s$ and sink $t$ as $G$ by Lemma~\ref{lem:samesink}. 
    We also have at least one route from $s$ to $t$, which come from the following two cases. 
    Either we have an internal vertex and so by Lemma~\ref{lem:innervertexroute} we have a route. 
    If we do not have an internal vertex, then we only have $s$ and $t$ and, as the graph is connected by Lemma~\ref{lem:equalsroute}, there must be at least one edge from $s$ to $t$. 
    So, let $R_2$ be a route from $s$ to $t$ using edges of $G-R_1$. 
    Continue constructing $R_i \in G - \{R_1, \ldots, R_{i-1}\}$ in this way. 
    Since removing a route strictly decreases the number of edges of the previous step, and the number of edges of $G$ is finite, this process must end. 
    The result is a decomposition $G = \cup_{i=1}^kR_i$ for some $k \in \mathbb{N}$ and moreover these routes are edge disjoint since they are edge disjoint at every step of the construction.
\end{proof}

\begin{corollary}\label{cor:numberofroutes}
    The number of routes in a route decomposition $\mathcal{R}$ is $\outdeg_G(s)$, where $s$ is the source of $G$.
\end{corollary}

\begin{proof}
    Let $k = \outdeg_G(s)$ and $\{e_i : i \in [k] \}$ be the edges out of $s$. Then since each edge of the graph is in some route from our decomposition, there exists a unique $R_i \in \mathcal{R}$ such that $e_i \in R_i$. 
    Similarly, given $R \in \mathcal{R}$, it must use exactly one edge $e$ out of $s$ and no other route can use $e$ as the routes have disjoint edge sets. 
    This correspondence shows $|\mathcal{R}| = k$.
\end{proof}

\begin{remark}\label{rem:multisources}
    We can also partition the edge set  of $G$ when $G$ has multiple sources and sinks by identifying all the sources and similarly identifying all the sinks. 
\end{remark}

\section{Equatorial Flow Triangulations}

Throughout this section, we assume that $G$ is a DAG which has no idle edges and satisfies degree equality, hence has a Gorenstein flow polytope.
First, note that the vector 
\[
\sum_{e\in E(G)}\chi_e=\mathbbm{1}=(1,1,\ldots,1)\in \R^{E(G)}
\]
is the Gorenstein point for the flow polytope of $G$, where $\chi_e$ denotes the indicator vector with a $1$ in the $e$-th position and $0$'s elsewhere.
This follows from the observation by von~Bell et al. ~\cite{vonbell2024triangulations} that if $\flow(G)$ is Gorenstein, then it has Gorenstein point 
\[
\mathbbm{1} \in \mathcal{F}(G)\cong \cone{\flow(G)}
\]
under this integral equivalence.

The following theorem provides a combinatorial method for producing special simplices in Gorenstein flow polytopes.

\begin{theorem}\label{thm:specialroutesimplex}
    If $G$ satisfies degree equality, then any route decomposition $\mathcal{R}$ corresponds to a special simplex in $\flow(G)$ and the codegree of $\flow(G)$ is $|\mathcal{R}|$.
\end{theorem}

\begin{proof}
     Given a route decomposition $\mathcal{R}$, the sum of indicator vectors of routes $\chi_R$ for $R \in \mathcal{R}$ gives
     \[
        \sum_{R \in \mathcal{R}} \chi_R = \sum_{e\in E(G)}\chi_e \, ,
     \]
     since our route decomposition partitions the edge set. 
     This is the Gorenstein point of $\flow(G)$ and so by Proposition~\ref{prop:specialsimplex}, we have that $\conv{\chi_R : R \in \mathcal{R}}$ is a special simplex.
     By Theorem~\ref{thm:codegreeinterior}, we obtain the claimed value of the codegree.
\end{proof}

We denote the special simplices corresponding to a route decomposition in the following manner.

\begin{definition}\label{def:routesimplex}
    Given a route decomposition $\mathcal{R}$, the corresponding special simplex
    \[
        \rsimplex = \conv{\chi_R: R\in \mathcal{R}}
    \]
    is called the \emph{route simplex} for $\mathcal{R}$.
\end{definition}

We now proceed to apply the Bruns-R\"omer theory from Section~\ref{sec:background} in the context of Gorenstein flow polytopes.
Our goal is to obtain a completely graph-theoretic description of the resulting equatorial flow triangulations.

\begin{definition}\label{def:equatorialcomplex}
    Suppose we have a route decomposition $\mathcal{R}$ for a DAG $G$ satisfying degree equality.
    Let $\Gamma$ be the polyhedral subcomplex of $\partial \flow(G)$ generated by the faces of the form $\bigcap_{R \in \mathcal{R}}F_R$
    where each $F_R$ is a facet of $\flow(G)$ such that $\chi_R \notin F_R$. We call $\Gamma$ the \emph{equatorial complex} of the flow polytope.
\end{definition}

Our next goal is to give a combinatorial criteria to determine when a set of routes are contained in a common face of the equatorial complex for $\flow(G)$.

\begin{definition}\label{def:transversal}
    Given a route decomposition $\mathcal{R}$ of $G$, a set of edges $M$ is a \emph{transversal of $\mathcal{R}$} if $M$ consists of exactly one edge from each route of $\mathcal{R}$.
    Equivalently, a transversal is a system of distinct representatives for the routes in $\mathcal{R}$. 
    We say a route in $G$ \emph{avoids the transversal $M$} if the route does not use the edges in $M$, i.e., the route carries zero flow on the edges in $M$. 
\end{definition}

\begin{figure}
    \centering
    \input{transversalexample.tikz}
    \caption{Given the route decomposition from Figure~\ref{fig:routedecomp}, the collection of dashed edges give a transversal $M$, and all routes of the DAG using only solid edges avoids $M$.}
    \label{fig:transversalexample}
\end{figure}
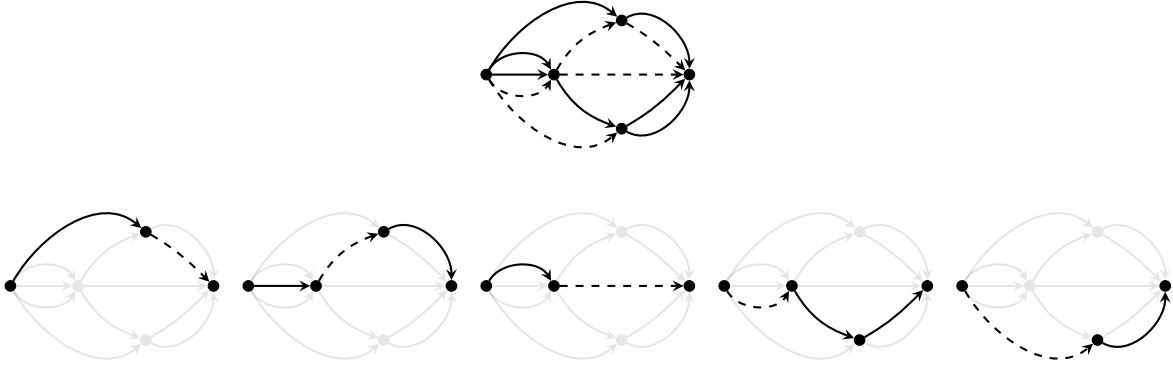

\begin{example}\label{ex:transversalexample}
    Figure~\ref{fig:transversalexample} shows an example of a transversal of a route decomposition.
\end{example}

The next few results give combinatorial criteria for describing both general faces and facets of Gorenstein flow polytopes.
Theorem~\ref{thm:facesavoidtransveral} gives our first combinatorial criterion for when a collection of routes form a common face of the equatorial complex for $\flow(G)$. 
In particular, we will see that transversals characterize the faces of $\Gamma$. 
Note that for a transversal $M=\{e_R : R \in \mathcal{R}\}$ we denote by $\bigcap_{R \in \mathcal{R}}\{x_{e_R} = 0\}$ the intersection of facets corresponding to $M$.

\begin{theorem}\label{thm:facesavoidtransveral}
    Given a directed acyclic graph $G$ with degree equality and a route decomposition $\mathcal{R}$, a collection of routes $\mathcal{S}$ form a face of $\Gamma$ if and only if $\mathcal{S}$ is the set of all routes avoiding a union of transversals in $\mathcal{R}$.
\end{theorem}
\begin{proof}
    We first will expand on Definition~\ref{def:equatorialcomplex}. 
    Facets of $\flow(G)$ are given by $x_e = 0$ for $e \in E(G)$. 
    Given a route $R \in \mathcal{R}$, since $\chi_R$ gives a flow of $1$ on the edges it uses in $G$, then $\chi_R$ is not on the facets given by $\{x_e = 0 : e \in E(R)\}$. 
    So, given $\{e_R \in E(R) : R \in \mathcal{R}\}$ to be any transversal, then
    by Definition~\ref{def:equatorialcomplex}, we have that
    \[
        \bigcap_{R \in \mathcal{R}}\{x_{e_R} = 0\}
    \]
    gives a generating face of the equatorial complex. 
    
    We will first show the backward direction. 
    Let $\mathcal{S}$ be the set of all routes avoiding a union of transversals. Let one such transversal be $M=\{e_R \in E(R) : R \in \mathcal{R}\}$, then the indicator vectors for routes avoiding $M$ are the vertices of the generating face of $\Gamma$ given by
    \[
        \bigcap_{R \in \mathcal{R}}\{x_{e_R} = 0\}
    \]
    since vertices of $\flow(G)$ and hence $\Gamma$ are routes.
    Similarly another transversal gives a generating face, and taking the intersection of all these generating faces is another face whose vertices are all routes avoiding the union of transversals.
    Thus $\mathcal{S}$ forms a face of $\Gamma$.
    For the forward direction, since our collection of routes form a face of the equatorial complex, the face with vertices $\mathcal{S}$ is the intersection of generating faces. Moreover each generating face is given by 
    \[
        \bigcap_{R \in \mathcal{R}} F_R
    \]
    where $F_R$ is a facet of $\flow(G)$ such that $\chi_R \notin F_R$. 
    We saw that facets of this form must use hyperplanes $\{x_e = 0 : e \in E(R)\}$. 
    As this holds for each route in our route decomposition, this gives a collection $\{e_R \in E(R) : R \in \mathcal{R}\}$ which is a transversal for this generating face. 
    Since routes in this generating face must give $0$ for these edges, these routes avoid this transversal. 
    As routes in $\mathcal{S}$ are all the vertices of the face formed by intersecting the generating faces, $\mathcal{S}$ is the set of all the routes which avoids the union of the transversals given by each generating face. 
\end{proof}

Corollary \ref{cor:facesavoidroutes} gives an alternative combinatorial characterization of when routes share a face of the equatorial complex. 

\begin{corollary}\label{cor:facesavoidroutes}
    Given a directed acyclic graph $G$ with degree equality, a collection of routes $\mathcal{S}$ lie on a common face of the equatorial complex if and only if there is no route in $\mathcal{R}$ contained in the union of edges of the routes in $\mathcal{S}$.
\end{corollary}
\begin{proof}
    We will prove the forward direction by contrapositive. 
    Suppose that some $R \in \mathcal{R}$ is the union of edges in $\mathcal{S}$, i.e. $E(R) \subseteq \bigcup_{S \in \mathcal{S}} E(S)$. 
    Then $\mathcal{S}$ must use all the edges of $R$ and so we cannot avoid a transversal of $\mathcal{R}$ since we cannot avoid any edge from $R$. 
    Thus, $\mathcal{S}$ does not lie on a common face by Theorem~\ref{thm:facesavoidtransveral}. 
    For the backward direction, now suppose no route from our route decomposition $\mathcal{R}$ in the union of edges of the routes in $\mathcal{S}$. 
    Then for every route $R \in \mathcal{R}$, there must exists some edge $e_R \in E(G)$ not used by $\mathcal{S}$, otherwise we would find a route as a union of edges. 
    This gives a transversal which routes in $\mathcal{S}$ avoid and so, by Theorem~\ref{thm:facesavoidtransveral}, these routes lie on a common face. 
\end{proof}

While Theorem~\ref{thm:facesavoidtransveral} describe exactly when a collection of routes form a face of the equatorial complex, not all these faces are facets. 
Theorem~\ref{thm:facettransversals} describes exactly when a transversal gives a facet of the equatorial complex.

\begin{theorem}\label{thm:facettransversals}
    Given a route decomposition of a directed acyclic graph $G$ and a generating face $F$ of $\Gamma$ avoiding a transversal $M=\{e_R : R \in \mathcal{R}\}$, $F$ is a facet of the equatorial complex $\Gamma$ if and only for all inner vertices $v$ of $G$ there exists a route incident to $v$ that does not use the edges in $M$. 
\end{theorem}
\begin{proof}
    We will prove the forward implication by contraspositive. 
    Suppose there exists an inner vertex $v \in V(G)$ such that all routes incident to $v$ use an edge in the transversal $M$. 
    Without loss of generality, suppose $R_1,\ldots, R_{\indeg_G(v)}$ are the routes in $\mathcal{R}$ using $v$.
    We will first show that the face of $\Gamma$ given by $\{x_{e_R} : R \in \mathcal{R}\}$ is contained in a face where all incoming edges of $v$ in $G$ are in the transversal, specifically, the face of $\Gamma$ given by $M'=\{e_R : R \in \mathcal{R} \setminus \{R_1,\ldots, R_{\indeg_G(v)}\} \} \cup \mathrm{in}(v)$. 
    Suppose we have a route $S$ which avoids the edges in  $M$. 
    Then our route $S$ would not be incident to $v$, as otherwise $S$ would use an edge of $M$. 
    This means that $E(S) \cap \mathrm{in}(v) = \emptyset$, thus $S$ also avoids the incoming edges of $v$, hence $S$ also avoids the edges in $M'$. 
    This shows our original face $F$ is contained in the modified face defined by the transversal $M'$.
    
    We will now show that the face of $\Gamma$ given by $M'$ is properly contained in another face and hence $F$ is not a facet. 
    Since $G$ has no idle edges, we have $\indeg_G(v) \geq 2$ and also $\outdeg_G(v) \geq 2$. 
    Moreover, for an incoming edge $e \in \mathrm{in}(v)$ there must be some outgoing edge $f \in \mathrm{out}(v)$ that lies on the same route, say $R_i$ for some $i$. Define a new transversal $M''=(M'\setminus \{e\}) \cup \{f\}$. 
    This contains all the routes avoiding the original transversal $M$, since every route avoiding $M$ is not incident to $v$. Moreover we claim this transversal contains a new route. 
    We have a path given by the route $R_i$ from $s$ to $v$; since $\outdeg_G(v) \geq 2$, there must be another route in $\mathcal{R}$, say $R_j$, using an edge out of $v$ different from $f$. 
    We can concatenate at $v$ the path from $s$ to $v$ following $R_i$ with the path from $v$ to $t$ following $R_j$ to create a new route that avoids $M''$. 
    Thus, $F$ is properly contained in another face of $\Gamma$ and hence is not a facet. 
    
    We will now prove the backward implication by contradiction. 
    For a contradiction, suppose $F$ avoiding the transversal $M$ is not a facet, and thus there exists a face $F'$ of $\Gamma$ such that $F \neq F'$ and $F \subset F'$. 
    Without loss of generality, $F'$ is a generating face of $\Gamma$ and so avoids a different transversal, say $\{e_R' : R \in \mathcal{R}\}$. Since we have a proper containment, there must exist a route $R' \in F'$ such that $R' \notin F$. 
    Since routes in $F$ avoid edges in $M$, there must exist an $S \in \mathcal{R}$ such that $e_{S} \in E(R')$ and $e_{S}' \notin E(R')$, since $R' \in F'$ for the same $S \in \mathcal{R}$. 
    This means in particular that $e_{S} \neq e_{S}'$. 
    We now claim there exists an inner vertex of indegree at least two in $V(S)$ appearing between $e_{S}$ and $e_{S}'$. 
    Since $e_{S} \neq e_{S}'$, there are at least three vertices and so the one between these edges in $S$ must be an inner vertex. 
    Since $G$ has no idle edges this inner vertex must have indegree greater than or equal to two.
    
    Having established our claim, suppose that $e_{S}$ appears before $e_{S}'$ in the route $S$ (when the order is switched, the following argument is symmetric). 
    Let $v$ be an inner vertex with $\indeg_G(v) = \outdeg_G(v) \geq 2$ between these two edges. 
    By our assumption, we have a route using $v$ in $G$ that avoids the edges in $M$. 
    In particular we have a path from the source to $v$ avoiding the edge $e_{S}$. 
    We also have a path from $v$ to the sink obtained by using the edges of $S$, which also avoids $e_{S}$ since this edge appears before $v$ in $S$. 
    Concatenating these paths we get a route from our source to sink which avoids $M$ and uses $e_{S}'$, hence the resulting route is contained in $F$ but not in $F'$. 
    This contradicts that $F \subset F'$, and hence $F$ is a facet.    
\end{proof}


We next introduce a triangulation of $\Gamma$ that uses as input only a route decomposition $\mathcal{R}$ of $G$ with a linear ordering of the routes in $\mathcal{R}$.

\begin{definition}
    Consider a route decomposition of $G$ with a linear ordering $\leq$ given by 
    \[
    \mathcal{R}=\{R_1<R_2<R_3<\cdots <R_k\}\, .
    \]
    For an inner vertex $v$ with indegree $j$ that is incident to the routes $R_{i_1}<R_{i_2}<\cdots <R_{i_j}$, with incoming edges at $v$ given by $e_{i_1},\ldots,e_{i_j}$, order the incoming edges according to the ordering of the routes, i.e., $e_{i_1}<\cdots <e_{i_j}$.
    We similarly order the outgoing edges at each vertex.
    This framing is called the \emph{route decomposition framing} and yields by Theorem~\ref{thm:dkk} a DKK triangulation $T$ of $\flow(G)$.
    The restriction $T \vert_\Gamma$ is denoted by $\dkkeq$.   
\end{definition}

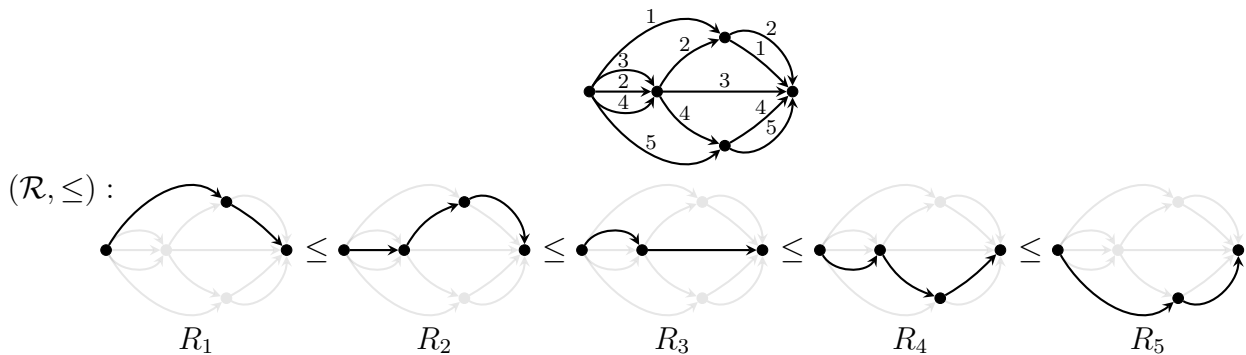
\begin{figure}
    \centering
    \input{routeframing.tikz}
    \caption{A route decomposition and framing.}
    \label{fig:routeframing}
\end{figure}

Note that it is often useful to consider the incoming edges at $v$ to be labeled $i_1,i_2,\ldots,i_j$, and similarly for outgoing edges.
A route decomposition and associated framing is given in Figure~\ref{fig:routeframing}.
Theorem~\ref{thm:gortriangulations} and Theorem~\ref{thm:specialroutesimplex} imply that the following triangulation exists and is well-defined.

\begin{definition}\label{def:equatorialflowtriangualtion}
    Given a route decomposition with a linear ordering $(\mathcal{R},\leq)$, we define the \emph{equatorial flow triangulation} to be the triangulation of $\flow(G)$ given by the join
    \[
        \dkkeq \ast \rsimplex \, ,
    \]
    i.e., the triangulation obtained as the union of simplices by $\conv{F \cup \rsimplex}$ where $F$ is a face of the triangulation $T_{\eq}$.
\end{definition}

The results of Bruns and R\"omer in Section~\ref{sec:background} yield the following in the special case of flow polytopes.

\begin{theorem}\label{thm:mainthm}
    For a Gorenstein flow polytope $\flow(G)$, the equatorial flow triangulation is regular and unimodular, and the equatorial sphere $T_{\eq}$ is the boundary complex of a simplicial polytope. 
    Further, the $h$-vector of $\dkkeq$ is the $h^*$-vector of $\flow(G)$.
\end{theorem}

\begin{proof}
This follows immediately from Theorem~\ref{thm:gortriangulations}, Theorem~\ref{thm:gorrestrictionsimplicial}, and Corollary~\ref{cor:gorensteinhhstar}.
\end{proof}

\begin{remark}
    We can replace $T_{\mathrm{eq}}$ by any regular unimodular triangulation of $\Gamma$, say $S$, to obtain a regular unimodular triangulation of $\flow(G)$ given by $S \ast \Delta_\mathcal{R}$.
    The usefulness and beauty of the equatorial flow triangulation is that it is defined entirely combinatorially by a route decomposition with a linear ordering.
\end{remark}

A purely combinatorial description of the equatorial flow triangulation is the following.

\begin{theorem}\label{thm:conditionsflowequatorial}
    For a directed acyclic graph $G$ satisfying degree equality and $(\mathcal{R},\leq)$ a route decomposition with a linear ordering, a collection of routes $C$ is a simplex in $\dkkeq$ if it satisfies precisely the following two conditions:
    \begin{enumerate}
        \item The routes in $C$ form a clique in the DKK triangulation induced by $\leq$ on $\mathcal{R}$.
        \item The routes in $C$ avoid a transversal $e_1,\ldots,e_k$ of $R_1, \ldots, R_k \in \mathcal{R}$.
    \end{enumerate}
\end{theorem}

\begin{proof}
    From the equatorial flow triangulation, $T_{eq}$ triangulates the equatorial complex $\Gamma$. 
    In particular, routes which form a simplex lie on a common face of $\Gamma$ and are a restriction of our DKK triangulation $T$ to $\Gamma$. 
    The first condition give routes which form a simplex from the DKK triangulation, while the second condition gives routes which lie on a common face by Theorem~\ref{thm:facesavoidtransveral}.
\end{proof}

What is particularly remarkable is that if $G$ satisfies degree equality and has a vertex of indegree $3$ or greater, then the equatorial flow triangulation is not a DKK triangulation, i.e., is not induced by a framing.


\begin{theorem}\label{thm:notfullnotdkk}
    If $G$ has an inner vertex with indegree at least $3$, then the equatorial flow triangulation $\dkkeq \ast \Delta_\mathcal{R}$ is not a DKK triangulation. 
\end{theorem}

\begin{proof}
    Suppose $v$ satisfies $\indeg_G(v) \geq 3$, hence by degree equality $\outdeg_G(v) \geq 3$. 
    This implies that there exist at least three routes in $\mathcal{R}$ that are incident to $v$, without loss of generality suppose we have at least the routes $R_1,R_2,R_3$.
    For the equatorial flow triangulation, $R_1, R_2, R_3$ are each adjacent to every other route, as they are in our route simplex. However, we will show that one of these three routes is in conflict with some other route, hence the equatorial flow triangulation is not a DKK triangulation.
    For $i \in[3]$, denote $e_{i}$ to be the incoming edge at $v$ and $f_{i}$ to be the outgoing edge at $v$ in $R_i$ as we have a route decomposition. 
    Since we have a framing, there exists permutations $\sigma, \pi \in S_3$ such that
    \[
        e_{\sigma(1)} < e_{\sigma(2)} < e_{\sigma(3)}, \quad f_{\pi(1)} < f_{\pi(2)} < f_{\pi(3)}.
    \]
    in our framing. There must exist a route $R_j$ for some $j \in [3]$ such that $e_j = e_{\sigma(2)}$. We will now consider the following two cases. 
    For the first case, suppose $f_j = f_{\pi(2)}$, then $R_j$ conflicts with the route given by the concatenation of the paths $R_{\sigma(1)}v$ and $v R_{\pi(3)}$.
    For the second case, suppose  $f_j \neq f_{\pi(2)}$ then either $f_j = f_{\pi(1)}$ or $f_j = f_{\pi(3)}$. 
    If $f_j = f_{\pi(1)}$, then $R_j$ conflicts with the route given by the concatenation of the paths $R_{\sigma(1)} v$ and $v R_{\pi(2)}$. 
    If $f_j = f_{\pi(3)}$, then $R_j$ conflicts with the route given by the concatenation of the paths $R_{\sigma(3)} v$ and $v R_{\pi(2)}$. 
    In any case, $R_j$ conflicts with a route, proving the claim. 
\end{proof}

Theorem \ref{thm:notfullnotdkk} does not cover DAGs satisfying degree equality where all inner vertices have indegree $2$. 
For such DAGs, the equatorial flow triangulation is not always the same triangulation as the DKK triangulation induced by the route decomposition framing, as shown by the following example that is equivalent to an example given by Reiner and Welker~\cite{reinerwelker} via the flow and order polytope correspondence given in Section~\ref{sec:order}.

\begin{figure}
    \centering
    \input{notDKK.tikz}
    \caption{A DAG with a route decomposition for which the equatorial flow triangulation and DKK triangulation are not the same.}
    \label{fig:notDKK}
\end{figure}
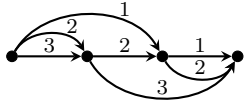

\begin{example}
    \label{ex:notDKK}
    For the DAG given in Figure~\ref{fig:notDKK}, in the equatorial flow triangulation, the route consisting of edges labeled by $2$ is connected to every other route by an edge.
    However, in the DKK triangulation for the route decomposition framing, the all-$2$'s route conflicts with the route along the spine labeled $321$, and hence do not form an edge.
\end{example}

\section{A Projection of the Flow Polytope}\label{sec:projection}

In this section, our goal is to describe an explicit geometric realization of the equatorial complex for a Gorenstein flow polytope $\flow(G)$ as the boundary of a reflexive polytope with the same $h^*$-polynomial as $\flow(G)$.
Specifically, we give both convex hull and (non-redundant) facet descriptions of the reflexive polytope. 
For an arbitrary polytope $P \subseteq \mathbb{R}^n$ and $V $ a linear subspace of $\mathbb{R}^n$, the \emph{quotient polytope} is
\[
    P/V := \{ p + V : p \in P\} \subseteq \mathbb{R}^n/V.
\]
In other words, the quotient polytope is the image of $P$ under the canonical quotient map $\mathbb{R}^n \to \mathbb{R}^n/V$. 
In addition, the quotient polytope in $\mathbb{R}^n/V$ is linearly isomorphic to the image of $P$ under any linear surjection $\mathbb{R}^n \to \mathbb{R}^{n - \dim V}$ with kernel $V$~\cite{athanasiadis}.
The following proposition is a mild restatement of results found in the work of Athanasiadis~\cite{athanasiadis} and Reiner and Welker~\cite{reinerwelker}.

\begin{prop}\label{prop:projection}
    Given a Gorenstein polytope $P \subseteq \mathbb{R}^n$ and a special simplex $\Delta$, let $V$ be the linear space spanned by $\Delta$. 
    The boundary complex of the quotient polytope $P/V$ is abstractly isomorphic to the equatorial complex $\Gamma(P,\Delta)$. 
    Moreover, if we have a triangulation of $P$ of the form $S \ast \Delta$, where $S$ is a triangulation of $\Gamma$, then the boundary complex of $P/V$ inherits a triangulation abstractly isomorphic to $S$.
\end{prop}

In particular, note that the facets of $P/V$ are in bijection with the facets of the equatorial complex.
Further, $P/V$ is known to be a reflexive polytope~\cite{brunsromer}.

We need to introduce some additional notation that will be used to describe the resulting reflexive polytope.
Assume we have a DAG $G$ (with no idle edges) satisfying degree equality with a route decomposition framing given by $R_1 < R_2 < \ldots <  R_k$. 
Given an edge $e \in E(G)$, let $\mathrm{label}(e)$ be the index of the route $R \in \mathcal{R}$ which contains this edge. 
This is well defined as we have a route decomposition. 
As usual, the inner vertices are labeled by $[n]$. 

\begin{definition}
    Given an internal vertex $i \in [n]$, let $\inlev(i)$ denote the set of labels of the incoming edges $e \in \mathrm{in}(i)$ at inner vertex $i$, i.e.
    \[
        \inlev(i) := \{ \mathrm{label}(e) : e \in \mathrm{in}(i) \}.
    \]
    Additionally, for a subset of edges $D \subseteq E(G)$, define
    \[
        \inlev(D,i) := \{ \mathrm{label}(e) : e \in \mathrm{in}(i) \cap D \}
    \]
    In particular, if we have a route $R$ define
    $\inlev(R,i) := \inlev(E(R),i)$. 
    Note that $\inlev(R,i)$ consists of at most a single element.
    
    Similarly, define $\outlev(i)$ as the set of labels from the outgoing edges $e \in \mathrm{out}_G(i)$, i.e.,
    \[
        \outlev(i) := \{ \mathrm{label}(e) : e \in \mathrm{out}(i)\}.
    \]
    We also similarly define $\outlev(R,i)$ for a route $R$ such that $i \in [n]$.
\end{definition}

\begin{definition}
    For each inner vertex $i \in [n]$ define
    \[
        V^i := \mathbb{R}^{\indeg(i)}
    \]
    to be the vector space with standard basis vectors $e^i_{j_1},\ldots,e^i_{j_\indeg(i)}$ where $\{j_1,\ldots,j_{\indeg(i)}\} = \inlev(i)$ and where $j_1 < j_2 < \ldots < j_{\indeg(i)}$ under the natural ordering of $\mathbb{N}$. Note that we use the superscript to denote that this basis is with respect to inner vertex $i$. 
\end{definition}

Given this notation, we are now able to define a linear map $\varphi$ with kernel given by the linear span of the special simplex formed  by the indicator vectors of $\{R_1,\ldots,R_k\}$.
The map given in Definition~\ref{def:phimap} is a generalization of a map defined by Jonah Berggren~\cite{berggrenpersonal} for flow polytopes of amply framed DAGs.

\begin{definition}\label{def:phimap}
    Define a linear map \[
    \varphi:\mathbb{R}^{|E(G)|} \to V^1 \times V^2 \times \cdots \times  V^n
    \]
    by extending linearly the map defined on a basis element $e_{(i,j)}$ for $(i,j)\in E(G)$ by
    \[
    \varphi(e_{(i,j)}):=e^j_{\mathrm{label}((i,j))}-e^i_{\mathrm{label}((i,j))}
    \]
    where if $i$ or $j$ is a source or a sink, we omit the corresponding summand.
\end{definition}

Observe that for a route $R$ in $G$, this definition of $\varphi$ yields
\begin{equation}\label{eqn:phiroute}        \varphi(R)=\sum_{i \in [n] \cap V(R)} e^i_{\inlev(R,i)} - e^i_{\outlev(R,i)} \, .
\end{equation}

\begin{definition}
    Given a DAG $G$ satisfying degree equality and an ordered route decomposition $\mathcal{R}$, we define 
    \[Q_\mathcal{R}:=\varphi(\flow(G)) \, .
    \]
\end{definition}

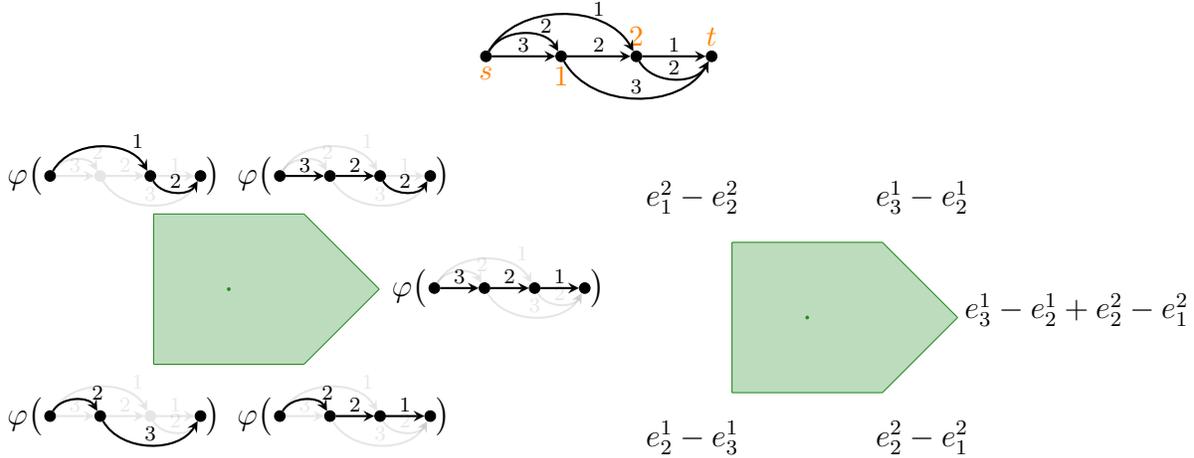
\begin{figure}
    \centering
    \input{quotientexample.tikz}
    \caption{A DAG with a route decomposition given by the edge labels, vertices of the DAG in orange, and the quotient polytope $Q_\mathcal{R}$.}
    \label{fig:quotientexample}
\end{figure}

\begin{figure}
    \centering
    \input{quotient3dexample.tikz}
    \caption{Examples of three dimensional quotient polytopes $Q_\mathcal{R}$ where the labels correspond to the route decompositions $\mathcal{R}$.}
    \label{fig:quotient3dexample}
\end{figure}
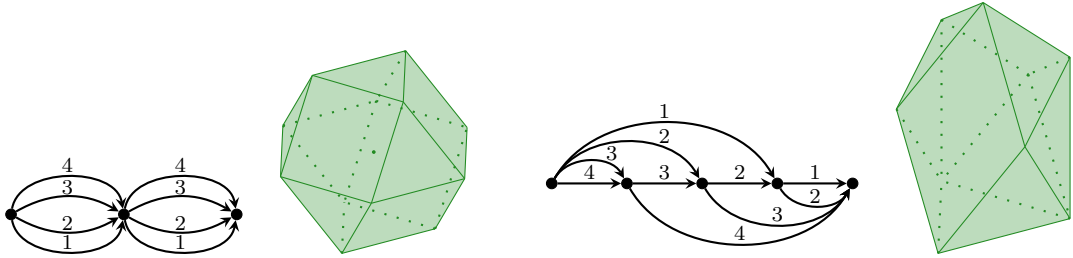

\begin{example}    Figure~\ref{fig:quotientexample} contains a DAG with a route decomposition and the corresponding polytope $Q_{\mathcal{R}}$ with vertices labeled by both routes $S$ and vectors $\varphi(S)$.
Figure~\ref{fig:quotient3dexample} provides two examples of DAGs with three-dimensional projection polytopes $Q_{\mathcal{R}}$.
Note that the left example in Figure~\ref{fig:quotient3dexample} generalizes to give the type-A root polytope as a quotient of a product of two simplices, i.e., suppose we have the standard basis vectors $e_k \in \mathbb{R}^n$ then the quotient polytope is given by $\conv{e_i - e_j : i \neq j \in [n]}$.
\end{example}

The following lemma connects the geometric structure of $Q_{\mathcal{R}}$ with the equatorial complex.

\begin{lemma}\label{lem:phi}
    The map $\varphi$ vanishes on the linear subspace spanned by the route simplex. 
    Thus, $Q_\mathcal{R}$ is a reflexive polytope whose face structure is given by the equatorial complex $\Gamma$ and whose vertices are 
    $\{\varphi(R):R\notin \mathcal{R}\}$.
\end{lemma}

\begin{proof}
That $\varphi$ vanishes on a route $R_j\in \mathcal{R}$ follows from the fact that all the edge labels in $R_j$ are $j$, and thus Equation~\eqref{eqn:phiroute} yields a telescoping sum that cancels to zero.
The claims regarding the face structure and vertices of $Q_{\mathcal{R}}$ follow immediately from Proposition~\ref{prop:projection}.
\end{proof}

Having given a convex-hull description of the reflexive projection $Q_{\mathcal{R}}$ of $\flow(G)$ for a special simplex defined by a route decomposition $\mathcal{R}$, our next task is to give a hyperplane description of $Q_{\mathcal{R}}$.
The following definition and theorem are the key tools we will use for this task.

\begin{definition}
    Given a transversal $\{e_R : R \in \mathcal{R}\}$, for every $R \in \mathcal{R}$ we have a partition of the edges $E(R)$ given by the following three sets
    \[
    R^{\mathrm{left}} \sqcup \{e_R\} \sqcup R^{\mathrm{right}}  = E(R)
    \]
    where $R^{\mathrm{left}}$ are all the edges in $G$ which the route $R$ uses before the edge $e_R$ and $R^{\mathrm{right}}$ are all the edges appearing in $R$ after the edge $e_R$. 
    Note that sometimes these sets can be empty.
\end{definition}

Given a route $S$ in our DAG $G$, we denote by $(\varphi(S))^i_l$ the $l$-th coordinate of $\varphi(S)$ in $V^i$.

\begin{theorem}\label{thm:trasnversalequation}
    For a DAG $G$ satisfying degree equality, a route decomposition $\mathcal{R}$ of $G$, and a transversal $M:=\{e_R : R \in \mathcal{R}\}$, every route $S$ satisfies  \begin{equation}\label{eq:transversalequation}
         \sum_{i \in [n]}\sum_{\substack{R \in \mathcal{R}, \\ l \in \inlev(R^{\mathrm{left}},i)}} (\varphi(S))^i_{l} = 1 - |E(S) \cap M| \, .
    \end{equation}
    In particular, for any choice of transversal $M$, every route $S$ satisfies the linear inequality
    \[
    \sum_{i \in [n]}\sum_{\substack{R \in \mathcal{R}, \\ l \in \inlev(R^{\mathrm{left}},i)}} x^i_{l} \, \leq 1 \, .
    \]
    Note that the edges in each $R^{\mathrm{left}}$ depend on the choice of $M$.
\end{theorem}

\begin{proof}
    Suppose $S \in \mathcal{R}$, then it must use exactly one edge in the transversal. Moreover, the image $\varphi(S) = 0$ and thus 
    \[
         \sum_{i \in [n]}\sum_{\substack{R \in \mathcal{R}, \\ l \in \inlev(R^{\mathrm{left}},i)}} (\varphi(S))^i_{l} = 0 = 1 - |E(S) \cap M| \, .
    \]
    Next, suppose $S \notin \mathcal{R}$. 
    Since we have a route decomposition, $S$ must exit the source following some $R \in \mathcal{R}$.
    Since $S\notin \mathcal{R}$, $S$ must eventually switch to follow edges in a different route in $\mathcal{R}$, and then possibly switch to follow another route in $\mathcal{R}$, and so forth.
    Hence, there is the following partition of $S$ as
    \[
        S = R_{i_1} v_1 R_{i_2} v_2 \cdots R_{i_{m-1}} v_{m-1} R_{i_m}
    \]
    for some integer $m$ where $m \geq 2$, where each $v_j \in [n]$ is an inner vertex where the edge into $v_j$ agrees with $R_{i_j} \in \mathcal{R}$ and the edge out of $v_j$ agrees with $R_{i_{j+1}} \in \mathcal{R}$ for $R_{i_j} \neq R_{i_{j+1}}$. 
    Since $S$ is not in our route simplex, at least one such $v_j$ exists. 
    
  This partition of the edges in $S$ implies that the image of $S$ under $\varphi$ is
    \[
        \varphi(S) = \sum_{j=1}^{m-1} e_{i_j}^{v_j} - e_{i_{j+1}}^{v_j} = (e^{v_1}_{i_1} - e^{v_1}_{i_2}) + (e^{v_2}_{i_2} - e^{v_2}_{i_3}) + \cdots + (e^{v_{m-1}}_{i_{m-1}} - e^{v_{m-1}}_{i_m})\, ,
    \]
    since all other vertices in $[n]$ are either not on $S$ and hence do not contribute to the sum or lie on $S$ but $S$ does not change levels at this vertex and hence there is no contribution for this vertex under the $\varphi$ map. 
    Because of this expression for $\varphi(S)$, we can simplify the left hand side of Equation~\eqref{eq:transversalequation} to only include the vertices $v_j$ for $j \in [m-1]$. This simplification is given by
    \begin{equation}\label{eq:transversalsimplified}
        \sum_{i \in [n]}  \sum_{\substack{R \in \mathcal{R}, \\ l \in \inlev(R^{\mathrm{left}},i)}} (\varphi(S))^i_{l} = \sum_{j \in [m-1]}  \sum_{\substack{R \in \mathcal{R}, \\ l \in \inlev(R^{\mathrm{left}},v_j)}} (\varphi(S))^{v_j}_{l} \, .
    \end{equation}
    Equation~\eqref{eq:transversalsimplified} also shows that we can compute all terms of this equation by considering each segment of $S$, i.e., on $R_{i_1} v_1$, on $v_{j}R_{i_{j+1}} v_{j+1}$ for all $j\in [m-1]$, and on $v_{m-1} R_{i_m}$, since the $j$-th term in the right-hand-side of Equation~\eqref{eq:transversalsimplified} is
    \[
        \sum_{\substack{R \in \mathcal{R}, \\ l \in \inlev(R^{\mathrm{left}},v_j)}} \hspace{-0.7cm}(\varphi(S))^{v_j}_{l} = 
        \begin{cases}
            (\varphi(S))_{i_j}^{v_j} & \text{if } i_j \in \inlev(R_{i_j}^{\mathrm{left}},v_j), \, i_{j+1} \notin \inlev(R_{i_{j+1}}^{\mathrm{left}},v_j),\\
            (\varphi(S))_{i_j}^{v_j} + (\varphi(S))_{i_{j+1}}^{v_j} & \text{if } i_j \in \inlev(R_{i_j}^{\mathrm{left}},v_j), \, i_{j+1} \in \inlev(R_{i_{j+1}}^{\mathrm{left}},v_j), \\
            (\varphi(S))_{i_{j+1}}^{v_j}  & \text{if } i_j \notin \inlev(R_{i_j}^{\mathrm{left}},v_j), \, i_{j+1} \in \inlev(R_{i_{j+1}}^{\mathrm{left}},v_j), \\
            0 & \text{otherwise.} 
        \end{cases} 
    \]
    For instance, when considering $v_1$, we can consider the first two segments $R_{i_1} v_1$ and $v_1 R_{i_2} v_2$ (or $R_{i_m}$ if there are only two segments in the partition). 
    Each of these make up the two possible values for the $j=1$ term in the sum above, and hence determining the values on each segment will compute Equation~\eqref{eq:transversalsimplified}. 

    We now use this setup to prove that Equation~\eqref{eq:transversalequation} holds.
    First consider $R_{i_1} v_1$. We have two cases. 
    In the first case, suppose $e_{R_{i_1}} \notin E(R_{i_1} v_1)$. 
    Then $E(R_{i_1} v_1) \subseteq R_{i_1}^{\mathrm{left}}$ since we start this path from the source. 
    Hence we have $(\varphi(S))_{i_1}^{v_1} = 1$ appearing in the right-hand sum for Equation~\eqref{eq:transversalsimplified}. 
    In the second case we have $e_{R_{i_1}} \in E(R_{i_1} v)$. 
    Then $i_1 \notin \inlev(R_{i_1}^{\mathrm{left}},v_1)$ since the edge into $v_1$ will either be $e_{R_{i_1}}$ or be in $R_{i_1}^{\mathrm{right}}$. 
    So we have a $0$ value, which we consider to be $1 - 1$. 

    We will now show that in the remaining segments, a $0$ is contributed to the right-hand side of Equation~\eqref{eq:transversalsimplified} if the segment contains no edge from the transversal and a $-1$ is contributed if an edge for the transversal is contained in the segment. 
    When combined with our analysis of the first segment, this will prove the equality in Equation~\eqref{eq:transversalequation}. 
    Consider $v_{j}R_{i_{j+1}} v_{j+1}$ for $j\in [m-1]$. 
    We again consider two cases. 
    For the first case, suppose $e_{R_{j+1}} \notin E(v_{j}R_{i_{j+1}} v_{j+1})$. 
    Then we consider two subcases.
    \begin{itemize}
        \item In the first subcase, if $E(v_j R_{i_{j+1}} v_{j+1}) \subseteq R^{\mathrm{left}}_{i_{j+1}}$, then $(\varphi(S))_{i_{j+1}}^{v_j} + (\varphi(S))_{i_{j+1}}^{v_{j+1}}$ 
        must appear in our $j$-th and $j+1$-st term of Equation~\eqref{eq:transversalsimplified}.
        However, this simplifies to
        \[
            (\varphi(S))_{i_{j+1}}^{v_j} + (\varphi(S))_{i_{j+1}}^{v_{j+1}} = (-e_{i_{j+1}}^{v_j})_{i_{j+1}}^{v_j} + (e_{i_{j+1}}^{v_{j+1}})_{i_{j+1}}^{v_{j+1}} = -1 + 1 = 0
        \]
        by using our segment expression of $\varphi(S)$. 

        \item In the second subcase, suppose $E(v_j R_{i_{j+1}} v_{j+1}) \not\subseteq R^{\mathrm{left}}_{i_{j+1}}$. 
        Since we know this segment does not contain $e_{R_{i_{j+1}}}$, then we must have $E(v_j R_{i_{j+1}} v_{j+1}) \subseteq R^{\mathrm{right}}_{i_{j+1}}$. 
        This means we have no terms associated to $(\varphi(S))_{i_{j+1}}^{v_j}$ or $(\varphi(S))_{i_{j+1}}^{v_{j+1}}$ and so there are no associated terms included in Equation~\eqref{eq:transversalsimplified}, contributing $0$.
    \end{itemize}
    Consider next the second case, where $e_{R_{j+1}} \in E(v_j R_{i_{j+1}} v_{j+1})$. 
    Then $i_{j+1} \in \inlev(R_{i_{j+1}}^{\mathrm{left}},v_j)$ and $i_{j+1} \notin \inlev(R_{i_{j+1}}^{\mathrm{left}},v_{j+1})$, since $e_{R_{i_{j+1}}}$ cannot appear in $R_{i_{j+1}}$ before $v_{j}$ since it is used in $v_{j}R_{i_{j+1}} v_{j+1}$ and the edge of $S$ into $v_{j+1}$ must be either $e_{R_{i_{j+1}}}$ or in $R_{i_{j+1}}^{\mathrm{right}}$. 
    This means we have 
    \[
        (\varphi(S))_{i_{j+1}}^{v_j} = (- e_{i_{j+1}}^{v_j}) = -1
    \]
    appearing for the $j$-th term of the right-hand side of Equation~\eqref{eq:transversalsimplified}.
    
    Our final step is to consider the last segment $v_{m-1} R_{i_m}$. 
    Again we have two cases. 
    In the first case, consider $e_{R_{i_m}} \notin E(v_{m-1} R_{i_m})$. 
    We must have $E(v_{m-1}R_{i_m}) \subseteq R^{\mathrm{right}}_{i_m}$ since $S$ avoids our transversal edge and this segment goes to the sink. 
    So, we do not get any additional contribution to the right-hand side of Equation~\eqref{eq:transversalsimplified}.
    In the second case, suppose $e_{R_{i_m}} \in E(v_{m-1} R_{i_m})$. 
    This implies that the edge of $R_{i_m}$ into $v_{m-1}$ must be in $R_{i_m}^{\mathrm{left}}$ and hence $i_m \in \inlev(R_{i_m}^{\mathrm{left}},v_{m-1})$. 
    So, we have the value
    \[
        (\varphi(S))_{i_m}^{v_{m-1}} = (- e_{i_m}^{v_{m-1}})_{i_m}^{v_{m-1}} = -1
    \]
    contributed to Equation~\eqref{eq:transversalsimplified} in the $v_{m-1}$-th summand. 

    Combining these observations shows that for the right-hand side of Equation~\eqref{eq:transversalsimplified}), we start with a value of $1$ coming from $R_{i_1} v_1$, and any time we have an edge from the transversal contained in a segment we add a value of $-1$ to this sum, while in the cases where we do not use a transversal edge, we do not change the value. 
    Hence, we have established that  
    \[
         \sum_{i \in [n]}\sum_{\substack{R \in \mathcal{R}, \\ l \in \inlev(R^{\mathrm{left}},i)}} (\varphi(S))^i_{l} = 1 - |E(S) \cap M|.
    \]
    
\end{proof}

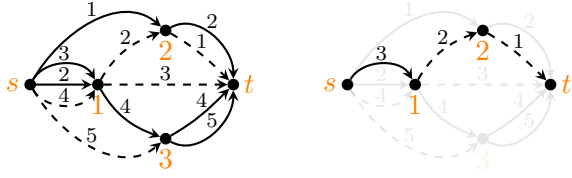
\begin{figure}
    \centering
    \input{routeequation.tikz}
    \caption{An example of a route and a DAG illustrating Equation~\eqref{eq:transversalequation}.}
    \label{fig:routeequation}
\end{figure}

\begin{example}\label{ex:routeequation}
    Consider the DAG with transversal and route decomposition given in Figure~\ref{fig:routeequation}.
    For the route $S$ following the edges labeled in the order $321$, we have
    \[
    \varphi(S)=e^1_3-e^1_2+e^2_2-e^2_1 \, ,
    \]
    and thus
    \[
    \sum_{i \in [n]}\sum_{\substack{R \in \mathcal{R}, \\ l \in \inlev(R^{\mathrm{left}},i)}} (\varphi(S))^i_{l}= (\varphi(S))^1_3 +(\varphi(S))^1_2 + (\varphi(S))^2_1 =1-1-1=-1=1-|E(S)\cap M|\, .
    \]
\end{example}

Combining our previous results, we are able to obtain the following halfspace description of $Q_{\mathcal{R}}$.

\begin{theorem}\label{thm:facetsofquotient}
For a DAG $G$ satisfying degree equality, with a route decomposition $\mathcal{R}$ of $G$, the polytope $Q_\mathcal{R}=\varphi(\flow(G))$ is full-dimensional in the subspace given by the intersection of
$\sum_{j=1}^{\indeg(i)}  x_j^i= 0$ for all $i \in [n]$, with facets given by the following halfspaces, one for each transversal $M:=\{e_R : R \in \mathcal{R}\}$ satisfying the conditions of Theorem~\eqref{thm:facettransversals}:
    \[
        H_M = \left\{ x \in V^1 \times \ldots \times V^n \: \middle\vert \: \sum_{i \in [n]}  \sum_{\substack{R \in \mathcal{R}, \\ l \in \inlev(R^{\mathrm{left}},i)}} x^i_{l}  = 1 \right\}
    \]
\end{theorem}

\begin{proof}
    Using Equation~\eqref{eqn:phiroute}, it follows that $\varphi(\flow(G))$ satisfies $\sum_{j=1}^{\indeg(i)}  x_j^i= 0$ for all $i \in [n]$.
    We will now show $\varphi(\flow(G))$ is full-dimensional in this space. 
    By Proposition~\ref{prop:projection} and Theorem~\ref{thm:gortriangulations}, the equatorial triangulation is given by a simplex joined with a triangulated sphere. 
    Suppose $m$ is the dimension of the equatorial complex (as a sphere) and $k$ is the out-degree of the source $k :=\outdeg(s)$. 
    As our equatorial flow triangulation is given as a join, we have
    \begin{align*}
        \dim(\flow(G)) &= \dim(\Delta_\mathcal{R}) + m + 1 \\
        &= (k-1) + m + 1 \\
        &= k + m \, .
    \end{align*}
    Since $\dim(\flow(G)) = |E(G)| - |\{ \text{inner vertices of } G\}| - 1$ then we can compute $m$ as
    \begin{align*}
        m &= |E(G)| - |\{ \text{inner vertices of } G\}| - 1 - k \\
        &= \left( \sum_{v \in V(G)} \outdeg(v) \right) - |\{ \text{inner vertices of } G\}| - 1 - k \\
        &= \left(\sum_{v \in \mathrm{inner} \, V(G)} \outdeg(v) \right)+ (\outdeg(s) - k) - |\{ \text{inner vertices of } G\}| - 1 \\
        &=\left(\sum_{v \in \mathrm{inner} \, V(G)} (\indeg(v) - 1) \right) - 1.
    \end{align*}
    The second equality holds since the each edge $(i,j) \in E(G)$ appears exactly once as the outgoing edge of the left endpoint $i$. 
    The third equality holds by removing the source term of the sum. 
    This term also gives $0$ when subtracting $k$. 
    The last equality holds by degree equality and since we subtract the number of inner vertices and so can subtract $1$ to each term of the sum. 
    Since $m$ is the dimension of the boundary of $Q_\mathcal{R}$ by proposition~\ref{prop:projection}, then 
    \[
        \dim Q_\mathcal{R} = m+1 = \sum_{v \in \mathrm{inner} \, V(G)} (\indeg(v) - 1)
    \]
    which is the dimension of the space in $V^1 \times \cdots \times V^n$ given by the intersection of the hyperplanes $\sum_{j=1}^{\indeg(i)}  x_j^i= 0$ for all $i \in [n]$. 
    Hence $Q_\mathcal{R}$ is full dimensional in this space.
    
    By Proposition~\ref{prop:projection}, the facets of $Q_{\mathcal{R}}$ correspond to the facets of the equatorial complex $\Gamma$.
    We know that each transversal satisfying the conditions of Theorem~\ref{thm:facettransversals} yields a facet of the equatorial complex.
    Suppose we have $H_M$ for some $M:=\{e_R : R \in \mathcal{R}\}$ satisfying the conditions of Theorem~\ref{thm:facettransversals}. 
    We claim any route $S$ avoids the transversal $M$ if and only if $\varphi(S) \in H_M$. 
    If $S$ avoids the transversal $M$, by Theorem~\ref{thm:trasnversalequation}, we have $\varphi(S) \in H_M$. 
    Similarly, if $\varphi(S) \in H_M$, then $S$ avoids $M$ since otherwise by Theorem~\ref{thm:trasnversalequation} the vector $\varphi(S)$ would not satisfy the defining equation of $H_M$. 
    Thus, a route $S$ is on the facet of $\Gamma$ corresponding to $M$ if and only if $\varphi(S)$ is on the hyperplane $H_M$, which by Theorem~\ref{thm:trasnversalequation} is a supporting hyperplane of $Q_{\mathcal{R}}$.
    This completes the proof.
\end{proof}

\section{Triangulations of Order Polytopes for Strongly Planar Posets}\label{sec:order}

In this section, we prove that for strongly planar DAGs with no idle edges, a particular equatorial flow triangulation is integrally equivalent to the equatorial triangulation of the associated order polytope defined by Reiner-Welker. 
Towards this goal, we first give the integral equivalence of flow polytopes and order polytopes when we have strongly planar posets and DAGs shown by M\'esz\'aros, Morales, and Striker.
Second, we give the definitions and results of Reiner-Welker to see the structure of the equatorial triangulation of order polytopes, before proving our result.

\begin{definition}\label{def:stronglyplanardag}
A DAG $G$ is \emph{strongly planar} if $G$ is a planar graph with a planar realization such that the following two conditions hold. First, for every directed edge $(i,j)$ of $G$, the $x$-coordinate of $i$ is strictly less than the $x$-coordinate of $j$.
Second, each edge $(i,j)$ is embedded in the plane as the graph of a piecewise differentiable function of $x$.
\end{definition}

\begin{definition}\label{def:stronglyplanarposet}
A poset $P$ is \emph{strongly planar} if the Hasse diagram of $P$ is a planar graph with a planar realization such that the following two conditions hold. First, for every directed edge $(i,j)$ in the Hasse diagram, the $y$-coordinate of $i$ is strictly less than the $y$-coordinate of $j$.
Second, each edge $(i,j)$ is embedded in the plane as the graph of a piecewise differentiable function of $y$.
\end{definition}

The following definition was originally given by M\'esz\'aros, Morales, and Striker~\cite[Section 3.3]{meszaros-morales-striker}.

\begin{definition}
    \label{def:dualgraph}
    Given a strongly planar DAG $G$, we define the \emph{truncated dual graph} $G^*$ to be a strongly planar realization of the dual graph of $G$ with the vertex corresponding to the infinite region deleted.
    The graph $G^*$ is thus the Hasse diagram of a strongly planar poset that we denote by $P_G$.
    Similarly, given a strongly planar poset $P$, let $H$ be the Hasse diagram of $P\cup \{\hat{0},\hat{1}\}$.
    We define $G_P$ to be the strongly planar DAG obtained as the strongly planar dual of $H$.
\end{definition}

Note that the definitions of planar DAG and strongly planar poset given in~\cite{meszaros-morales-striker} are not sufficient for Theorem~\ref{thm:flowandorderequivalence} to hold in general. 
Definitions~\ref{def:stronglyplanardag} and~\ref{def:stronglyplanarposet} require additional structure for planar realizations of edges that are sufficient to establish Theorem~\ref{thm:flowandorderequivalence} using the proofs given in~\cite{meszaros-morales-striker}.

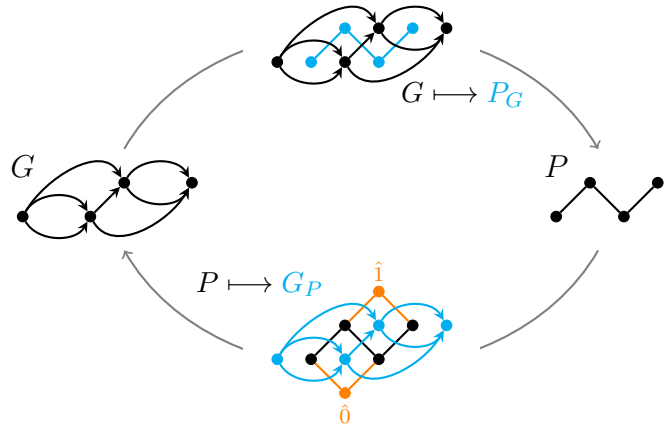
\begin{figure}
    \centering
    \input{stronglyplanar.tikz}
    \caption{A strongly planar poset $P$ and a strongly planar DAG $G$ which are related by duality.}
    \label{fig:stronglyplanar}
\end{figure}

\begin{example}
    \label{ex:stronglyplanar}
    Figure~\ref{fig:stronglyplanar} gives an example of a strongly planar DAG $G$ and poset $P$. Moreover, these are related by duality and the corresponding constructions of $G^*$ and $H$ are given. 
\end{example}

We next recall basic facts about order polytopes of posets.

\begin{definition}\label{def:orderpolytope}
Given a finite poset $P = \{t_1,\ldots,t_n\}$, define the \emph{order polytope} of $P$ to be 
\[
\op{P}:= \{(x_1,\ldots,x_n) \in [0,1]^{|P|}:x_i\leq x_j \text{ if }t_i\leq_P t_j\} \, .
\]
\end{definition}

\begin{remark}\label{rem:ordercone}
    Given a poset $P$, a map $f:P \to \mathbb{R}_{\geq 0}$ is \emph{order preserving} if $f(t_i) \leq f(t_j)$ whenever $t_i \leq_P t_j$. 
    The set of order preserving maps for $P$ forms a polyhedral cone called the \emph{order cone}, denoted $K_P$. 
    This means the order polytope can be seen as bounding the codomain of maps in $K_p$ to be in $[0,1]$ since we can set $f(t_i) := x_i$ and similarly for $x_j$.
\end{remark}

\begin{definition}
    An \emph{upper order ideal} (or \emph{filter}) $F$ is a subset of a poset $P$ with the property that if $x \in F$ and $x \leq_P y$ for some $y \in P$, then $y \in F$. Let $J(P)$ denote the poset of upper order ideals ordered by inclusion.
\end{definition}

The following proposition gives the vertex description of an order polytope.

\begin{prop}\label{prop:verticesorderpolytope}
    Given a finite poset $P$, the vertices of the order polytope $\mathcal{O}(P)$ are characteristic functions of upper order ideals, i.e., for  $F \subseteq P$ an upper order ideal the order preserving map $\chi_F:P \to \mathbb{R}_{\geq 0}$ defined by
    \[
        \chi_F(p) = 
        \begin{cases}
            1 & \text{if } p \in F, \\
            0 & \text{otherwise}
        \end{cases}
    \]
    gives a vertex of $\mathcal{O}(P)$.
\end{prop}

The following theorem is the key correspondence between flow and order polytopes in the strongly planar case.

\begin{theorem}[M\'esz\'aros, Morales, and Striker~\cite{meszaros-morales-striker}]\label{thm:flowandorderequivalence}
Given a strongly planar DAG $G$, the flow polytope $\flow(G)$ is integrally equivalent to the order polytope of the strongly planar poset dual to $G$.
Conversely, given a strongly planar poset $P$, the order polytope $\mathcal{O}(P)$ is integrally equivalent to the flow polytope of the strongly planar DAG $G_P$.
\end{theorem}

As a consequence of Theorems~\ref{thm:mainthm} and~\ref{thm:flowandorderequivalence}, we have the following.

\begin{lemma}\label{lem:routedecompisgraded}
    A strongly planar DAG $G$ has a route decomposition if and only if $P_G$ is graded.
\end{lemma}
\begin{proof}
    A strongly planar DAG $G$ has a route decomposition if and only if $\flow(G)$ is Gorenstein. By Theorem~\ref{thm:flowandorderequivalence} this holds if and only if $\Or(P_G)$ is Gorenstein which holds if and only if $P_G$ is a graded poset. 
\end{proof}

 The following theorem shows how our equatorial flow triangulations can be used to produce new equatorial triangulations of Gorenstein order polytopes for strongly planar posets, describable in purely graph theoretic language.
 
\begin{theorem}
    If $G$ is a strongly planar DAG satisfying degree equality, then every route decomposition of $G$ induces a unimodular triangulation of the corresponding order polytope $\mathcal{O}(P_G)$.
\end{theorem}

In addition to the integral equivalence given in Proposition~\ref{prop:canonicaldkk}, M\'esz\'aros, Morales, and Striker give an equivalence of the canonical triangulation of order polytopes and a DKK triangulation given by a framing called the planar framing. 
The structure of the canonical triangulation is nicely given by the following definition and theorem. 

\begin{definition}\label{def:ordercomplex}
    The order complex of $J(P)$ is defined as the simplicial complex $\Delta(J(P))$ where each upper order ideal $F$ of $P$ defines a vertex and each chain $F_1 \subseteq \cdots \subseteq F_t$ of nested ideals defines a simplex.
\end{definition}

\begin{theorem}\label{thm:cannonicalchains}
    The canonical triangulation of the order polytope is isomorphic as an abstract simplicial complex to the order complex of $J(P)$, i.e., chains of nested upper order ideals correspond to the simplicies of the canonical triangulation. 
\end{theorem}

 We will now define the framing which gives an equivalence of the DKK triangulation and the canonical triangulation for a strongly planar DAG. 

\begin{definition}[M\'esz\'aros, Morales, and Striker~\cite{meszaros-morales-striker}]\label{def:planarframing}
    Consider a strongly planar DAG $G$, then for every inner vertex there is an order of the edges given by the planar embedding. At inner vertex $j$ we can order the incoming and outgoing edges from top to bottom. More specifically, draw a small circle around the inner vertex such that each edge hits the circle once at a point and order the incoming and outgoing edges by decreasing $y$-coordinates of these points. This framing is called the \emph{planar framing} of $G$ (as it depends on the planar embedding of the DAG $G$).
\end{definition}

\begin{prop}[M\'esz\'aros, Morales, and Striker~\cite{meszaros-morales-striker}]\label{prop:canonicaldkk}
    Consider a strongly planar DAG $G$. Under the integral equivalence map $\phi$ given by M\'esz\'aros, Morales, and Striker, the canonical triangulation of $\mathcal{O}(P_G)$ maps to a DKK triangulation of $\flow(G)$ given by the planar framing 
\end{prop}

A key observation for the proof of M\'esz\'aros, Morales, and Striker is given by the following proposition. It will also be useful later in this section in proving Theorem~\ref{thm:equatorialisequatorialflow}.

\begin{prop}[M\'esz\'aros, Morales, and Striker~\cite{meszaros-morales-striker}]\label{prop:idealsandroutes}   
    Suppose we have a strongly planar DAG $G$ and let $\phi$ be the integral equivalence of the order polytope and flow polytope given by M\'esz\'aros, Morales, and Striker.
    Given a vertex $\chi_F$ of $\mathcal{O}(P_G)$ for an upper order ideal $F$ of $P_G$, we have $\phi(\chi_F)$ is the unit flow along a route in $G$ separating $F$ and $P_G - F$. Moreover, any route $R$ in $G$ separates some upper order ideal $F$ and $P_G - F$.
    
\end{prop}

We now discuss the equatorial triangulation defined by Reiner-Welker~\cite{reinerwelker} on order polytopes. The definitions used in this paper are order preserving, while Reiner and Welker used order reversing. However, this does not change their construction of the equatorial triangulation for the order polytope. We give the corresponding definitions here.

Let $P$ be a graded poset with $r$ rank sets, starting with rank $1$. Let $P_i$ denote the set of elements of $P$ of rank $i$ for $0 \leq i \leq r$ where we say $P_0 = \emptyset$.

\begin{definition}\label{def:rankconstant}
    An order preserving map $f:P \to \mathbb{R}_{\geq 0}$ is \emph{rank-constant} if it is constant along ranks, in other words $f(p) = f(q)$ whenever $p,q \in P_i$ for some $i \in [r]$
\end{definition}

\begin{definition}\label{def:equatorialideal}
    An order preserving map $f:P \to \mathbb{R}_{\geq 0}$ is \emph{equatorial} if $\min_{p \in P} f(p) = 0$ and for every $j \in [2,r]$ there exists a cover relation $p_{j-1} \prec_P p_j$ such that $p_{j-1} \in P_{j-1}$, $p_j \in P_j$, and $f(p_{j-1}) = f(p_j)$.
\end{definition}

\begin{definition}\label{def:rankorequatorialchains}
    Suppose we have a chain of upper order ideals $F_1 \subset \cdots \subset F_t$. We say this is a \emph{rank-constant chain} (resp. \emph{equatorial chain}) if the sum $\chi_{F_1} + \cdots + \chi_{F_t}$ is rank-constant (resp. equatorial)
\end{definition}

\begin{example}\label{ex:rankconstantideals}
    The only rank constant upper order ideals are of the form
    \[
        F_j^{rc} = \bigsqcup_{i > j} P_i
    \]
    for $0 \leq j \leq r$ where $F_r^{rc} = \emptyset$.
\end{example}

\begin{definition}\label{def:jumps}
    Given a chain of ideals $F_1 \subset \cdots \subset F_t$, we say its \emph{jumps} are 
    \[
        J_i = F_i - F_{i-1}
    \]
    for $i \in [t+1]$ and where $F_0 =\emptyset$, $F_{t+1} = P$ 
\end{definition}

\begin{theorem}[Reiner and Welker~\cite{reinerwelker}]\label{thm:equatorialchains}
    A chain of nonempty upper order ideals $F_1 \subset \cdots \subset F_t$ is equatorial if and only if its jumps have the following property; for every $j \in [2,r]$ there exists $p_{j-1} \prec_P p_j$ with $p_{j-1} \in P_{j-1}$ $p_j \in P_j$ and a value $i \in [t+1]$ such that $p_{j-1}, p_j \in J_i$. In other words, between every pair of consecutive ranks there exists a cover relation which appears in the same jump.
\end{theorem}

\begin{definition}\label{def:equatorialsubcomplex}
    We define $\Delta_{eq}(P)$ as the subcomplex of the order complex $\Delta(J(P))$ whose faces are given by equatorial chains of non-empty ideals.
\end{definition}

\begin{remark}\label{rem:equatorialcovnention}
    Reiner and Welker \cite{reinerwelker} define $\Delta_{eq}(P)$ as the equatorial complex which is different from our convention in defining $\Gamma$ to be the equatorial complex. We will see later that $\Delta_{eq}(P)$ is in fact related to $T_{eq}$.
\end{remark}

Reiner and Welker define the equatorial triangulation and show it has a nice structure given by a join. 

\begin{corollary}[Reiner and Welker~\cite{reinerwelker}]\label{cor:posetequatorialjoin}
    The equatorial triangulation of the order polytope $\mathcal{O}(P)$ is abstractly isomorphic to the simplicial join of $\sigma \ast \Delta_{eq}(P)$ where $\sigma$ is the interior simplex spanned by the rank constant ideas $\{F_j^{rc}\}_{j=0}^r$.
\end{corollary}

This join structure  of the equatorial triangulation of the order polytope is very similar to the equatorial flow triangulations. 
In fact, Theorem~\ref{thm:equatorialisequatorialflow} shows that the equatorial triangulation of the order polytope is equivalent to an equatorial flow triangulation, under the integral equivalence of flow polytopes and order polytopes for a strongly planar DAG.

\begin{remark}\label{rem:codomainofintegralequivalence}
Define the poset $\widehat{P_G}:= P_G \cup \{\hat{0}, \hat{1} \}$. 
    The integral equivalences of $\flow(G)$ and $\mathcal{O}(P_G)$ given by M\'esz\'aros, Morales, and Striker \cite{meszaros-morales-striker} are defined as follows:
    \begin{itemize}
        \item Define $\varphi:\flow(G) \to \Or(P_G)$ by $fl \mapsto (f(x))_{x \in P_G}$ where $f:P_G \to \mathbb{R}_{\geq 0}$ is given by
        \[
            f(x) = \sum_{e} fl(e)
        \]
        where we sum over the edges $e \in G$ which correspond to edges from a chain from $\hat{0}$ to $x$ in $\widehat{P_G}$ by duality. This works for any chain.

        \item Define $\phi:\Or(P_G) \to \flow(G)$ by $(f(x))_{x \in P_G} \to fl$ where $fl:E(G) \to \mathbb{R}_{\geq 0}$ is given by
        \[
            fl(e) = f(y) - f(x)
        \]
        where the edge $e \in G$ corresponds to the cover relation $x \prec y$ in $\widehat{P_G}$ by duality and $f(\hat{1}) = 1$, $f(\hat{0}) = 0$.
    \end{itemize}
    Note that the map the map $\varphi$ given by M\'esz\'aros, Morales, and Striker give in their paper is slightly stronger. In fact, the codomain of $\varphi$ is all order preserving functions $f$ of $P_G \cup \{\hat{0}, \hat{1} \}$ such that $f(\hat{0}) = 0$ and $f(\hat{1}) = 1$.
\end{remark}

\begin{theorem}\label{thm:equatorialisequatorialflow}
    Suppose $G$ is strongly planar and $\flow(G)$ is Gorenstein, then the equatorial triangulation of $\Or(P_G)$ corresponds to an equatorial flow triangulation under the integral equivalence of $\flow(G)$ and $\Or(P_G)$.
\end{theorem}
\begin{proof}
    Define the poset $\widehat{P_G}:= P_G \cup \{\hat{0}, \hat{1} \}$. 
    We will first show the rank constant simplex for the equatorial triangulation of $\Or(P_G)$ is equivalent to the route simplex for a particular choice of route decomposition of $\flow(G)$. 
    We have $P_G$ is graded by Lemma~\ref{lem:routedecompisgraded}, say of rank $r$. Recall that $(P_G)_k$ denotes the set of elements of $P_G$ of rank $k$ for $1 \leq k \leq r$. 
    Since $G$ is strongly planar, each edge of our DAG is a piecewise differentiable function. 
    So for each $x$-value between the $x$-values of $s$ and $t$, there exists a largest $y$-value attained by some edge in $G$. 
    By planarity, the collection of edges achieving these maximum values forms a route $R_0$ where all other edges in $G$ have smaller $y$-values. 
    By duality, the edges of $R_0$ correspond to the cover relations of $\widehat{P_G}$ between $\hat{1}$ and elements of $(P_G)_r$. 
    Moreover, by Proposition~\ref{prop:idealsandroutes}, this gives the upper order ideal $\emptyset = F_r^{rc}$ of $P_G$. 
    We can now consider $G - R_0$ and repeat the argument. 
    This gives a route $R_1$ where all the edges correspond to the cover relations of elements between $(P_G)_r$ and $(P_G)_{r-1}$. 
    Hence by Proposition~\ref{prop:idealsandroutes}, this gives the ideal $(P_G)_r = F_{r-1}^{rc}$ of $P_G$. 
    We continue in this way, and the last route to remove will be the route $R_r$ whose edges correspond to the cover relations between $(P_G)_1$ and $\hat{0}$ in $\widehat{P_G}$ which gives the ideal $P_G = F_{0}^{rc}$ of $P_G$. 
    So our route decomposition and linear order $(\mathcal{R}, \leq)$ is given by $R_0 \leq R_1 \leq \cdots \leq R_r$ which corresponds to the simplex given by rank constant ideals $F_r^{rc} \subset F_{r-1}^{rc} \subset \cdots \subset F_0^{rc}$. 
    In addition, the route decomposition framing given by $\mathcal{R}$ is equivalent to the planar framing as we will be ordering edges from the top down given by the planar embedding. 
    
    Let $T$ be the route decomposition framing given by $(\mathcal{R}, \leq)$, i.e., the planar framing. 
    We will now show that a simplex is in $T_{eq}$ of the flow polytope $\flow(G)$ if and only if the corresponding simplex under the integral equivalence to $\mathcal{O}(P_G)$ is in $\Delta_{eq}(P)$. 
    A simplex in $T_{eq}$ is a clique of routes under the planar framing which must also avoid a transversal by Theorem~\ref{thm:conditionsflowequatorial} and equivalence to the route decomposition framing. 
    By Proposition~\ref{prop:canonicaldkk}, we have a clique of routes if and only if the corresponding vertices in order polytope $\mathcal{O}(P_G)$ give a chain of nested upper order ideals. 
    Thus, we can assume coherency and it suffices to check that routes avoid a transversal if and only if for the corresponding ideals in $P_G$ satisfy the equatorial condition and are non-empty to prove our main claim. 
        
    For the forward direction, suppose we have coherent routes $\{S_i\}_{i=1}^t$ which avoid a transversal $M = \{e_R \in E(G): R \in \mathcal{R}\}$. 
    Consider the ideals $\{\varphi(S_i)\}_{i=1}^t$ of $P_G$ and let $f_{R_i}$ denote the cover relation $b_i \prec a_i$ of $\widehat{P_G}$ corresponding to the edge $e_{R_i} \in M$. 
    We will first show these ideals are nonempty. 
    Consider $e_{R_0}$, we have a cover relation $b_0 \prec a_0 = \hat{1}$ in $\widehat{P_G}$, which gives
    \[
        \varphi(S_i)(\hat{1}) - \varphi(S_i)(b_0) = \chi_{S_i}(e_{R_0}) + \varphi(S_i)(b_0) - \varphi(S_i)(b_0) = 0
    \]
    by definition of $\varphi$ for all $i$. 
    Note in this case we are using the extended codomain of order preserving functions discussed in Remark~\ref{rem:codomainofintegralequivalence}. We have shown $1 = \varphi(S_i)(\hat{1}) =  \varphi(S_i)(b_0)$. So $b_0 \in P_G$ is an element of all ideals $\{\varphi(S_i)\}_{i=1}^t$ and hence they are nonempty.
    We will now verify the equatorial condition. So for all $j \in [r-1]$, each $b_i\prec a_j$ is also a relation in $P_G$. 
    Hence for all $i$ we find the following property of $\varphi$. 
    When considering the edge of the cover relation $b_j \prec a_j$, we have
    \[
        \varphi(S_i)(a_j) - \varphi(S_i)(b_j) = \chi_{S_i}(e_{R_i}) + \varphi(S_i)(b_j) - \varphi(S_i)(b_j) = 0.
    \]
    Hence, $\varphi(S_i)(a_j) = \varphi(S_i)(b_j)$. 
    Since the edges of $R_j \in \mathcal{R}$ corresponds with the cover relations between consecutive ranks in $P_G$ by prior observations, then we have cover relations $f_{R_j}$ for all $j \in [r-1]$ between every pair of consecutive ranks of $P_G$. 
    Moreover, when we consider the edge $e_{R_r} \in M$ this corresponds to a cover relation $b_r=\hat{0} \prec a_r$ for some $a_r \in P_G$. Hence 
    \[
        \varphi(S_i)(a_r)= \chi_{S_i}(e_{R_r}) = 0
    \]
    meaning the minimum value of all $\varphi(S_i)$ is $0$ and achieved at $a_r \in P_G$. 
    Hence, $\sum_{i=1}^t \chi_{\varphi(S_i)}$
    has minimum value $0$ at $a_r$ and is constant on $f_{R_i}$ for $i \in [r-1]$ which shows $\{\varphi(S_i)\}_{i=1}^t$ is an equatorial chain of nonempty ideals.

    For the backward direction, suppose we have an equatorial chain of nonempty upper order ideals $F_1 \subset \cdots \subset F_t$ of $P_G$. 
    Since $F_1$ is nonempty and an upper order ideal, there exists $q_r \in F_1$ of the top rank. 
    Define $e_{R_0}$ to be the edge in $G$ corresponding to the cover relation $q_r \prec \hat{1}$ of $\widehat{P_G}$. 
    Moreover since we have a chain of ideals, $q_r \in \bigcap_{i=1}^t F_i$ and
    \[
        \phi(F_i)(e_{R_0}) = 1 - \chi_{F_i}(q_r) = 0
    \]
    for all $i \in [t]$. 
    So the routes $\{\phi(F_i)\}_{i=1}^t$ avoid the edge $e_{R_0}$ which is an edge in $R_0$. 
    Next, the equatorial condition gives $p_{j-1} \prec p_j$ in $P_G$ for $j \in[2,r]$ where $p_k \in (P_G)_k$ such that $\chi_{F_1}+ \cdots + \chi_{F_t}$ is constant. 
    By reindexing, we have for $k \in [r-1]$ the cover relations $p_{r-k} \prec p_{r-k+1}$ in $P_G$ such that $\chi_{F_1}+ \cdots + \chi_{F_t}$ is constant. 
    Define the edge $e_{R_k}$ in $G$ to be the one dual to $p_{r-k} \prec p_{r-k+1}$ of $P_G$. 
    By our construction of the route decomposition, we saw that the edges in $G$ corresponding to cover relations between $(P_G)_{r-k}$ and $(P_G)_{r-k+1}$ are the edges of $R_k \in \mathcal{R}$.
    Theorem~\ref{thm:equatorialchains} implies that $p_{r-k}$ and $p_{r-k+1}$ must either both be elements of $F_i$ or both not be elements of $F_i$ for all $i \in[t]$.  
    By definition of $\phi$, the flow on the edge $e_{R_k}$ for any $k \in [r-1]$ is
    \[
        \phi(F_i)(e_{R_k}) = \chi_{F_i}(p_{r-k+1}) - \chi_{F_i}(p_{r-k}) = 0
    \]
    by our previous observation of $p_{r-k}$ and $p_{r-k+1}$. To define the edge $e_{R_r}$, since we have an equatorial chain we know the minimum value of $\chi_{F_1}+ \cdots + \chi_{F_t}$ is $0$. So there must exist an element in $P_G$ of rank $1$ that attains this value, otherwise this would imply the poset only has nonzero values as we have upper order ideals. Let this element be $q_1 \in P_G$ hence $q_1$ is not in any of our ideals. Define the edge $e_{R_r} \in G$ be the one dual to the relation $\hat{0} \prec q_1$, so then 
    \[
        \phi(F_i)(e_{R_0}) = \chi_{F_i}(q_1) - 0 = 0.
    \]
    Thus, the routes $\{\phi(F_i)\}_{i=1}^t$ avoid the transversal given by $\{e_{R_j}\}_{j=0}^r$. 

    Therefore, we have shown $T_{eq}$ is integrally equivalent to $\Delta_{eq}$. 
    This establishes an isomorphism of triangulations $\sigma \ast \Delta_{eq}(P)$ and $\Delta_\mathcal{R} \ast T_{eq}$, since both are joins of equivalent simplicial complexes.
\end{proof}

For strongly planar graded posets, the equatorial flow triangulation is given by one route decomposition, but any route decomposition gives an equatorial flow triangulation. 
Hence, it is important to keep in mind that using equatorial flow triangulations we can find many different equatorial triangulations of order polytopes in the strongly planar case, and describe them combinatorially using DAGs.

\bibliographystyle{plain}
\bibliography{refs}

\end{document}

%% file: clique.tikz
\begin{tikzpicture}
    \begin{scope}[xshift=0, yshift=0, scale=1]    
        \vertex[fill](a1) at (1,0) {};
        \vertex[fill](a2) at (2,0) {};
        \vertex[fill](a3) at (3,0) {};
        \vertex[fill](a4) at (4,0) {};

        
        \draw[-stealth, thick] (a1) to[out=60,in=120]  (a3) node[xshift=-3, yshift=11] {\tiny $1$};
        \draw[-stealth, thick] (a1) to[out=60,in=120]  (a2) node[xshift=-5, yshift=11] {\tiny $1$};
        \draw[-stealth, thick] (a1) to[out=0,in=180] (a2) node[xshift=-10, yshift=4] {\tiny $2$} ;

        \draw[-stealth, thick] (a2) node[xshift=6, yshift=4] {\tiny $1$} to[out=0,in=180] (a3) node[xshift=-10, yshift=4] {\tiny $2$} ;
        \draw[-stealth, thick] (a2) node[xshift=8, yshift=-4.5] {\tiny $2$} to[out=-60,in=-120] (a4);

        \draw[-stealth, thick] (a3) node[xshift=6, yshift=4] {\tiny $1$} to[out=0,in=180] (a4);
        \draw[-stealth, thick] (a3) node[xshift=10, yshift=-4] {\tiny $2$} to[out=-60,in=-120]  (a4);
    \end{scope}
    \begin{scope}[xshift=-100, yshift=-50, scale=1]    
        \vertex[fill](a1) at (1,0) {};
        \vertex[fill](a2) at (2,0) {};
        \vertex[fill](a3) at (3,0) {};
        \vertex[fill](a4) at (4,0) {};

        
        \draw[-stealth, thick, opacity=0.1] (a1) to[out=60,in=120]  (a3) node[xshift=-3, yshift=11] {\tiny $1$};
        \draw[-stealth, thick] (a1) to[out=60,in=120]  (a2) node[xshift=-5, yshift=11] {\tiny $1$};
        \draw[-stealth, thick, opacity=0.1] (a1) to[out=0,in=180] (a2) node[xshift=-10, yshift=4] {\tiny $2$} ;

        \draw[-stealth, thick] (a2) node[xshift=6, yshift=4] {\tiny $1$} to[out=0,in=180] (a3) node[xshift=-10, yshift=4] {\tiny $2$} ;
        \draw[-stealth, thick, opacity=0.1] (a2) node[xshift=8, yshift=-4.5] {\tiny $2$} to[out=-60,in=-120] (a4);

        \draw[-stealth, thick] (a3) node[xshift=6, yshift=4] {\tiny $1$} to[out=0,in=180] (a4);
        \draw[-stealth, thick, opacity=0.1] (a3) node[xshift=10, yshift=-4] {\tiny $2$} to[out=-60,in=-120]  (a4);
    \end{scope}
    \begin{scope}[xshift=0, yshift=-50, scale=1]    
        \vertex[fill](a1) at (1,0) {};
        \vertex[fill](a2) at (2,0) {};
        \vertex[fill](a3) at (3,0) {};
        \vertex[fill](a4) at (4,0) {};

        
        \draw[-stealth, thick, opacity=0.1] (a1) to[out=60,in=120]  (a3) node[xshift=-3, yshift=11] {\tiny $1$};
        \draw[-stealth, thick, opacity=0.1] (a1) to[out=60,in=120]  (a2) node[xshift=-5, yshift=11] {\tiny $1$};
        \draw[-stealth, thick] (a1) to[out=0,in=180] (a2) node[xshift=-10, yshift=4] {\tiny $2$} ;

        \draw[-stealth, thick] (a2) node[xshift=6, yshift=4] {\tiny $1$} to[out=0,in=180] (a3) node[xshift=-10, yshift=4] {\tiny $2$} ;
        \draw[-stealth, thick, opacity=0.1] (a2) node[xshift=8, yshift=-4.5] {\tiny $2$} to[out=-60,in=-120] (a4);

        \draw[-stealth, thick] (a3) node[xshift=6, yshift=4] {\tiny $1$} to[out=0,in=180] (a4);
        \draw[-stealth, thick, opacity=0.1] (a3) node[xshift=10, yshift=-4] {\tiny $2$} to[out=-60,in=-120]  (a4);
    \end{scope}
    \begin{scope}[xshift=100, yshift=-50, scale=1]    
        \vertex[fill](a1) at (1,0) {};
        \vertex[fill](a2) at (2,0) {};
        \vertex[fill](a3) at (3,0) {};
        \vertex[fill](a4) at (4,0) {};

        
        \draw[-stealth, thick, opacity=0.1] (a1) to[out=60,in=120]  (a3) node[xshift=-3, yshift=11] {\tiny $1$};
        \draw[-stealth, thick, opacity=0.1] (a1) to[out=60,in=120]  (a2) node[xshift=-5, yshift=11] {\tiny $1$};
        \draw[-stealth, thick] (a1) to[out=0,in=180] (a2) node[xshift=-10, yshift=4] {\tiny $2$} ;

        \draw[-stealth, thick] (a2) node[xshift=6, yshift=4] {\tiny $1$} to[out=0,in=180] (a3) node[xshift=-10, yshift=4] {\tiny $2$} ;
        \draw[-stealth, thick, opacity=0.1] (a2) node[xshift=8, yshift=-4.5] {\tiny $2$} to[out=-60,in=-120] (a4);

        \draw[-stealth, thick, opacity=0.1] (a3) node[xshift=6, yshift=4] {\tiny $1$} to[out=0,in=180] (a4);
        \draw[-stealth, thick] (a3) node[xshift=10, yshift=-4] {\tiny $2$} to[out=-60,in=-120]  (a4);
    \end{scope}
    \begin{scope}[xshift=-50, yshift=-90, scale=1]    
        \vertex[fill](a1) at (1,0) {};
        \vertex[fill](a2) at (2,0) {};
        \vertex[fill, opacity=0.1](a3) at (3,0) {};
        \vertex[fill](a4) at (4,0) {};

        
        \draw[-stealth, thick, opacity=0.1] (a1) to[out=60,in=120]  (a3) node[xshift=-3, yshift=11] {\tiny $1$};
        \draw[-stealth, thick, opacity=0.1] (a1) to[out=60,in=120]  (a2) node[xshift=-5, yshift=11] {\tiny $1$};
        \draw[-stealth, thick] (a1) to[out=0,in=180] (a2) node[xshift=-10, yshift=4] {\tiny $2$} ;

        \draw[-stealth, thick, opacity=0.1] (a2) node[xshift=6, yshift=4] {\tiny $1$} to[out=0,in=180] (a3) node[xshift=-10, yshift=4] {\tiny $2$} ;
        \draw[-stealth, thick] (a2) node[xshift=8, yshift=-4.5] {\tiny $2$} to[out=-60,in=-120] (a4);

        \draw[-stealth, thick, opacity=0.1] (a3) node[xshift=6, yshift=4] {\tiny $1$} to[out=0,in=180] (a4);
        \draw[-stealth, thick, opacity=0.1] (a3) node[xshift=10, yshift=-4] {\tiny $2$} to[out=-60,in=-120]  (a4);
    \end{scope}
    \begin{scope}[xshift=50, yshift=-90, scale=1]    
        \vertex[fill](a1) at (1,0) {};
        \vertex[fill, opacity=0.1](a2) at (2,0) {};
        \vertex[fill](a3) at (3,0) {};
        \vertex[fill](a4) at (4,0) {};

        
        \draw[-stealth, thick] (a1) to[out=60,in=120]  (a3) node[xshift=-3, yshift=11] {\tiny $1$};
        \draw[-stealth, thick, opacity=0.1] (a1) to[out=60,in=120]  (a2) node[xshift=-5, yshift=11] {\tiny $1$};
        \draw[-stealth, thick, opacity=0.1] (a1) to[out=0,in=180] (a2) node[xshift=-10, yshift=4] {\tiny $2$} ;

        \draw[-stealth, thick, opacity=0.1] (a2) node[xshift=6, yshift=4] {\tiny $1$} to[out=0,in=180] (a3) node[xshift=-10, yshift=4] {\tiny $2$} ;
        \draw[-stealth, thick, opacity=0.1] (a2) node[xshift=8, yshift=-4.5] {\tiny $2$} to[out=-60,in=-120] (a4);

        \draw[-stealth, thick] (a3) node[xshift=6, yshift=4] {\tiny $1$} to[out=0,in=180] (a4);
        \draw[-stealth, thick, opacity=0.1] (a3) node[xshift=10, yshift=-4] {\tiny $2$} to[out=-60,in=-120]  (a4);
    \end{scope}
\end{tikzpicture}

%% file: degreeequalityDAG.tikz
\begin{tikzpicture}
     \begin{scope}[xshift=0, yshift=0, scale=1]
        \vertex[fill](a1) at (1,0) {};
        \vertex[fill](a2) at (2,0) {};
        \vertex[fill](a3) at (3,0.8) {};
        \vertex[fill](a4) at (3,-0.8) {};
        \vertex[fill](a5) at (4,0) {};
        
        \draw[-stealth, thick] (a1) to[out=0,in=180] (a2);
        \draw[-stealth, thick] (a1) to[out=60,in=120] (a2);
        \draw[-stealth, thick] (a1) to[out=-60,in=-120] (a2);
        \draw[-stealth, thick] (a1) to[out=60,in=140] (a3);
        \draw[-stealth, thick] (a1) to[out=-60,in=-140] (a4);
        
        \draw[-stealth, thick] (a2) to[out=60,in=200] (a3);
        \draw[-stealth, thick] (a2) to[out=-60,in=-200] (a4);
        \draw[-stealth, thick] (a2) to[out=0,in=180] (a5);
        
        \draw[-stealth, thick] (a3) to[out=30,in=90] (a5);
        \draw[-stealth, thick] (a3) to[out=-30,in=135] (a5);
        
        \draw[-stealth, thick] (a4) to[out=-30,in=-90] (a5);
        \draw[-stealth, thick] (a4) to[out=30,in=-135] (a5);
    
    \end{scope}     
        
\end{tikzpicture}

%% file: gorensteinface.tikz
\begin{tikzpicture}
     \begin{scope}[xshift=0, yshift=0, scale=1]
        \vertex[fill](a1) at (1,0) {};
        \vertex[fill](a2) at (2,0) {};
        \vertex[fill](a3) at (3,0) {};
        \vertex[fill](a4) at (4,0) {};
        
        \draw[-stealth, thick] (a1) to[out=0,in=-180] (a2);
        \draw[-stealth, thick] (a1) to[out=60,in=120] (a2);

        \draw[-stealth, thick] (a2) to[out=0,in=-180] (a3);

        \draw[-stealth, thick] (a3) to[out=0,in=-180] (a4);
        \draw[-stealth, thick] (a3) to[out=60,in=120] (a4);

        \node at (4.75,0) {$\longrightarrow$};
    \end{scope}  
    \begin{scope}[xshift=130, yshift=0, scale=1]
        \vertex[fill](a1) at (1,0) {};
        \vertex[fill](a2) at (2,0) {};
        \vertex[fill](a3) at (3,0) {};
        \vertex[fill](a4) at (4,0) {};
        
        \draw[-stealth, thick] (a1) to[out=0,in=-180] (a2);
        \draw[-stealth, thick] (a1) to[out=60,in=120] (a2);
        \draw[-stealth, thick,orange] (a1) to[out=60,in=120] (a3);

        \draw[-stealth, thick] (a2) to[out=0,in=-180] (a3);
        \draw[-stealth, thick,orange] (a2) to[out=-60,in=-120] (a4);

        \draw[-stealth, thick] (a3) to[out=0,in=-180] (a4);
        \draw[-stealth, thick] (a3) to[out=60,in=120] (a4);
    
    \end{scope} 
\end{tikzpicture}

%% file: routedecomp.tikz
\begin{tikzpicture}
     \begin{scope}[xshift=0, yshift=0, scale=0.9]
        \vertex[fill](a1) at (1,0) {};
        \vertex[fill](a2) at (2,0) {};
        \vertex[fill](a3) at (3,0.8) {};
        \vertex[fill](a4) at (3,-0.8) {};
        \vertex[fill](a5) at (4,0) {};
        
        \draw[-stealth, thick] (a1) to[out=0,in=180] (a2);
        \draw[-stealth, thick] (a1) to[out=60,in=120] (a2);
        \draw[-stealth, thick] (a1) to[out=-60,in=-120] (a2);
        \draw[-stealth, thick] (a1) to[out=60,in=140] (a3);
        \draw[-stealth, thick] (a1) to[out=-60,in=-140] (a4);
        
        \draw[-stealth, thick] (a2) to[out=60,in=200] (a3);
        \draw[-stealth, thick] (a2) to[out=-60,in=-200] (a4);
        \draw[-stealth, thick] (a2) to[out=0,in=180] (a5);
        
        \draw[-stealth, thick] (a3) to[out=30,in=90] (a5);
        \draw[-stealth, thick] (a3) to[out=-30,in=135] (a5);
        
        \draw[-stealth, thick] (a4) to[out=-30,in=-90] (a5);
        \draw[-stealth, thick] (a4) to[out=30,in=-135] (a5);
    
    \end{scope}   
    \begin{scope}[xshift=-180, yshift=-80, scale=0.9]
        \vertex[fill](a1) at (1,0) {};
        \vertex[fill,opacity=0.1](a2) at (2,0) {};%
        \vertex[fill](a3) at (3,0.8) {};
        \vertex[fill,opacity=0.1](a4) at (3,-0.8) {};%
        \vertex[fill](a5) at (4,0) {};
        
        \draw[-stealth, thick,opacity=0.1] (a1) to[out=0,in=180] (a2);%
        \draw[-stealth, thick,opacity=0.1] (a1) to[out=60,in=120] (a2);%
        \draw[-stealth, thick,opacity=0.1] (a1) to[out=-60,in=-120] (a2);%
        \draw[-stealth, thick] (a1) to[out=60,in=140] (a3);
        \draw[-stealth, thick,opacity=0.1] (a1) to[out=-60,in=-140] (a4);%
        
        \draw[-stealth, thick,opacity=0.1] (a2) to[out=60,in=200] (a3);%
        \draw[-stealth, thick,opacity=0.1] (a2) to[out=-60,in=-200] (a4);%
        \draw[-stealth, thick,opacity=0.1] (a2) to[out=0,in=180] (a5);%
        
        \draw[-stealth, thick,opacity=0.1] (a3) to[out=30,in=90] (a5);%
        \draw[-stealth, thick] (a3) to[out=-30,in=135] (a5);
        
        \draw[-stealth, thick,opacity=0.1] (a4) to[out=-30,in=-90] (a5);%
        \draw[-stealth, thick,opacity=0.1] (a4) to[out=30,in=-135] (a5);%
    \end{scope}  
    \begin{scope}[xshift=-90, yshift=-80, scale=0.9]
        \vertex[fill](a1) at (1,0) {};
        \vertex[fill](a2) at (2,0) {};
        \vertex[fill](a3) at (3,0.8) {};
        \vertex[fill, opacity=0.1](a4) at (3,-0.8) {};
        \vertex[fill](a5) at (4,0) {};
        
        \draw[-stealth, thick] (a1) to[out=0,in=180] (a2);
        \draw[-stealth, thick, opacity=0.1] (a1) to[out=60,in=120] (a2);%
        \draw[-stealth, thick, opacity=0.1] (a1) to[out=-60,in=-120] (a2);%
        \draw[-stealth, thick, opacity=0.1] (a1) to[out=60,in=140] (a3);%
        \draw[-stealth, thick, opacity=0.1] (a1) to[out=-60,in=-140] (a4);%
        
        \draw[-stealth, thick] (a2) to[out=60,in=200] (a3);
        \draw[-stealth, thick, opacity=0.1] (a2) to[out=-60,in=-200] (a4);
        \draw[-stealth, thick, opacity=0.1] (a2) to[out=0,in=180] (a5);
        
        \draw[-stealth, thick] (a3) to[out=30,in=90] (a5);
        \draw[-stealth, thick,opacity=0.1] (a3) to[out=-30,in=135] (a5);
        
        \draw[-stealth, thick, opacity=0.1] (a4) to[out=-30,in=-90] (a5);
        \draw[-stealth, thick, opacity=0.1] (a4) to[out=30,in=-135] (a5);
    \end{scope} 
    \begin{scope}[xshift=0, yshift=-80, scale=0.9]
        \vertex[fill](a1) at (1,0) {};
        \vertex[fill](a2) at (2,0) {};
        \vertex[fill, opacity=0.1](a3) at (3,0.8) {};
        \vertex[fill, opacity=0.1](a4) at (3,-0.8) {};
        \vertex[fill](a5) at (4,0) {};
        
        \draw[-stealth, thick, opacity=0.1] (a1) to[out=0,in=180] (a2);
        \draw[-stealth, thick] (a1) to[out=60,in=120] (a2);
        \draw[-stealth, thick, opacity=0.1] (a1) to[out=-60,in=-120] (a2);
        \draw[-stealth, thick, opacity=0.1] (a1) to[out=60,in=140] (a3);
        \draw[-stealth, thick, opacity=0.1] (a1) to[out=-60,in=-140] (a4);
        
        \draw[-stealth, thick, opacity=0.1] (a2) to[out=60,in=200] (a3);
        \draw[-stealth, thick, opacity=0.1] (a2) to[out=-60,in=-200] (a4);
        \draw[-stealth, thick] (a2) to[out=0,in=180] (a5);
        
        \draw[-stealth, thick, opacity=0.1] (a3) to[out=30,in=90] (a5);
        \draw[-stealth, thick, opacity=0.1] (a3) to[out=-30,in=135] (a5);
        
        \draw[-stealth, thick, opacity=0.1] (a4) to[out=-30,in=-90] (a5);
        \draw[-stealth, thick, opacity=0.1] (a4) to[out=30,in=-135] (a5);
    \end{scope} 
    \begin{scope}[xshift=90, yshift=-80, scale=0.9]
        \vertex[fill](a1) at (1,0) {};
        \vertex[fill](a2) at (2,0) {};
        \vertex[fill, opacity=0.1](a3) at (3,0.8) {};
        \vertex[fill](a4) at (3,-0.8) {};
        \vertex[fill](a5) at (4,0) {};
        
        \draw[-stealth, thick, opacity=0.1] (a1) to[out=0,in=180] (a2);
        \draw[-stealth, thick, opacity=0.1] (a1) to[out=60,in=120] (a2);
        \draw[-stealth, thick] (a1) to[out=-60,in=-120] (a2);
        \draw[-stealth, thick, opacity=0.1] (a1) to[out=60,in=140] (a3);
        \draw[-stealth, thick, opacity=0.1] (a1) to[out=-60,in=-140] (a4);
        
        \draw[-stealth, thick, opacity=0.1] (a2) to[out=60,in=200] (a3);
        \draw[-stealth, thick] (a2) to[out=-60,in=-200] (a4);
        \draw[-stealth, thick, opacity=0.1] (a2) to[out=0,in=180] (a5);
        
        \draw[-stealth, thick, opacity=0.1] (a3) to[out=30,in=90] (a5);
        \draw[-stealth, thick, opacity=0.1] (a3) to[out=-30,in=135] (a5);
        
        \draw[-stealth, thick,opacity=0.1] (a4) to[out=-30,in=-90] (a5);
        \draw[-stealth, thick] (a4) to[out=30,in=-135] (a5);
    \end{scope} 
    \begin{scope}[xshift=180, yshift=-80, scale=0.9]
        \vertex[fill](a1) at (1,0) {};
        \vertex[fill, opacity=0.1](a2) at (2,0) {};
        \vertex[fill, opacity=0.1](a3) at (3,0.8) {};
        \vertex[fill](a4) at (3,-0.8) {};
        \vertex[fill](a5) at (4,0) {};
        
        \draw[-stealth, thick, opacity=0.1] (a1) to[out=0,in=180] (a2);
        \draw[-stealth, thick, opacity=0.1] (a1) to[out=60,in=120] (a2);
        \draw[-stealth, thick, opacity=0.1] (a1) to[out=-60,in=-120] (a2);
        \draw[-stealth, thick, opacity=0.1] (a1) to[out=60,in=140] (a3);
        \draw[-stealth, thick] (a1) to[out=-60,in=-140] (a4);
        
        \draw[-stealth, thick, opacity=0.1] (a2) to[out=60,in=200] (a3);
        \draw[-stealth, thick, opacity=0.1] (a2) to[out=-60,in=-200] (a4);
        \draw[-stealth, thick, opacity=0.1] (a2) to[out=0,in=180] (a5);
        
        \draw[-stealth, thick, opacity=0.1] (a3) to[out=30,in=90] (a5);
        \draw[-stealth, thick, opacity=0.1] (a3) to[out=-30,in=135] (a5);
        
        \draw[-stealth, thick] (a4) to[out=-30,in=-90] (a5);
        \draw[-stealth, thick, opacity=0.1] (a4) to[out=30,in=-135] (a5);
    \end{scope}

\end{tikzpicture}

%% file: deleteroute.tikz
\begin{tikzpicture}
     \begin{scope}[xshift=0, yshift=0, scale=1]
        \vertex[fill](a1) at (1,0) {};
        \vertex[fill](a2) at (2,0) {};
        \vertex[fill](a3) at (3,0) {};
        \vertex[fill, color=orange](a4) at (1.5,0.3) {};
        
        \draw[-stealth, thick, color=orange] (a1) to[out=60,in=180] (a4);
        \draw[-stealth, thick] (a1) to[out=-60,in=-120] (a2);

        \draw[-stealth, thick, color=orange] (a4) to[out=0,in=120] (a2);

        \draw[-stealth, thick, color=orange] (a2) to[out=60,in=120] (a3);
        \draw[-stealth, thick] (a2) to[out=-60,in=-120] (a3);

        \node at (2,-0.75) {$G$};
        \node[color=orange] at (1.5,0.55) {$R_1$};
        \node at (3.75,0) {$\longrightarrow$};
    \end{scope}    
    \begin{scope}[xshift=100, yshift=0, scale=1]
        \vertex[fill](a1) at (1,0) {};
        \vertex[fill, color=cyan](a2) at (2,0) {};
        \vertex[fill](a3) at (3,0) {};
        
        \draw[-stealth, thick, color=cyan] (a1) to[out=-60,in=-120] (a2);

        \draw[-stealth, thick, color=cyan] (a2) to[out=-60,in=-120] (a3);

        \node at (2,-0.75) {$G - R_1$};
        \node[color=cyan] at (1.5,0) {$R_2$};
        \node at (3.75,0) {$\longrightarrow$};
        
    \end{scope}
    \begin{scope}[xshift=200, yshift=0, scale=1]
        \vertex[fill](a1) at (1,0) {};
        \vertex[fill](a3) at (3,0) {};
        

        \node at (2,-0.75) {$G - \{R_1, R_2\}$};
    \end{scope}
        
\end{tikzpicture}

%% file: transversalexample.tikz
\begin{tikzpicture}
     \begin{scope}[xshift=0, yshift=0, scale=0.9]
        \vertex[fill](a1) at (1,0) {};
        \vertex[fill](a2) at (2,0) {};
        \vertex[fill](a3) at (3,0.8) {};
        \vertex[fill](a4) at (3,-0.8) {};
        \vertex[fill](a5) at (4,0) {};
        
        \draw[-stealth, thick] (a1) to[out=0,in=180] (a2);
        \draw[-stealth, thick] (a1) to[out=60,in=120] (a2);
        \draw[-stealth, thick, dashed] (a1) to[out=-60,in=-120] (a2);
        \draw[-stealth, thick] (a1) to[out=60,in=140] (a3);
        \draw[-stealth, thick, dashed] (a1) to[out=-60,in=-140] (a4);
        
        \draw[-stealth, thick, dashed] (a2) to[out=60,in=200] (a3);
        \draw[-stealth, thick] (a2) to[out=-60,in=-200] (a4);
        \draw[-stealth, thick, dashed] (a2) to[out=0,in=180] (a5);
        
        \draw[-stealth, thick] (a3) to[out=30,in=90] (a5);
        \draw[-stealth, thick, dashed] (a3) to[out=-30,in=135] (a5);
        
        \draw[-stealth, thick] (a4) to[out=-30,in=-90] (a5);
        \draw[-stealth, thick] (a4) to[out=30,in=-135] (a5);
    
    \end{scope}   
    \begin{scope}[xshift=-180, yshift=-80, scale=0.9]
        \vertex[fill](a1) at (1,0) {};
        \vertex[fill,opacity=0.1](a2) at (2,0) {};%
        \vertex[fill](a3) at (3,0.8) {};
        \vertex[fill,opacity=0.1](a4) at (3,-0.8) {};%
        \vertex[fill](a5) at (4,0) {};
        
        \draw[-stealth, thick,opacity=0.1] (a1) to[out=0,in=180] (a2);%
        \draw[-stealth, thick,opacity=0.1] (a1) to[out=60,in=120] (a2);%
        \draw[-stealth, thick,opacity=0.1] (a1) to[out=-60,in=-120] (a2);%
        \draw[-stealth, thick] (a1) to[out=60,in=140] (a3);
        \draw[-stealth, thick,opacity=0.1] (a1) to[out=-60,in=-140] (a4);%
        
        \draw[-stealth, thick,opacity=0.1] (a2) to[out=60,in=200] (a3);%
        \draw[-stealth, thick,opacity=0.1] (a2) to[out=-60,in=-200] (a4);%
        \draw[-stealth, thick,opacity=0.1] (a2) to[out=0,in=180] (a5);%
        
        \draw[-stealth, thick,opacity=0.1] (a3) to[out=30,in=90] (a5);%
        \draw[-stealth, thick, dashed] (a3) to[out=-30,in=135] (a5);
        
        \draw[-stealth, thick,opacity=0.1] (a4) to[out=-30,in=-90] (a5);%
        \draw[-stealth, thick,opacity=0.1] (a4) to[out=30,in=-135] (a5);%
    \end{scope}  
    \begin{scope}[xshift=-90, yshift=-80, scale=0.9]
        \vertex[fill](a1) at (1,0) {};
        \vertex[fill](a2) at (2,0) {};
        \vertex[fill](a3) at (3,0.8) {};
        \vertex[fill, opacity=0.1](a4) at (3,-0.8) {};
        \vertex[fill](a5) at (4,0) {};
        
        \draw[-stealth, thick] (a1) to[out=0,in=180] (a2);
        \draw[-stealth, thick, opacity=0.1] (a1) to[out=60,in=120] (a2);%
        \draw[-stealth, thick, opacity=0.1] (a1) to[out=-60,in=-120] (a2);%
        \draw[-stealth, thick, opacity=0.1] (a1) to[out=60,in=140] (a3);%
        \draw[-stealth, thick, opacity=0.1] (a1) to[out=-60,in=-140] (a4);%
        
        \draw[-stealth, thick, dashed] (a2) to[out=60,in=200] (a3);
        \draw[-stealth, thick, opacity=0.1] (a2) to[out=-60,in=-200] (a4);
        \draw[-stealth, thick, opacity=0.1] (a2) to[out=0,in=180] (a5);
        
        \draw[-stealth, thick] (a3) to[out=30,in=90] (a5);
        \draw[-stealth, thick,opacity=0.1] (a3) to[out=-30,in=135] (a5);
        
        \draw[-stealth, thick, opacity=0.1] (a4) to[out=-30,in=-90] (a5);
        \draw[-stealth, thick, opacity=0.1] (a4) to[out=30,in=-135] (a5);
    \end{scope} 
    \begin{scope}[xshift=0, yshift=-80, scale=0.9]
        \vertex[fill](a1) at (1,0) {};
        \vertex[fill](a2) at (2,0) {};
        \vertex[fill, opacity=0.1](a3) at (3,0.8) {};
        \vertex[fill, opacity=0.1](a4) at (3,-0.8) {};
        \vertex[fill](a5) at (4,0) {};
        
        \draw[-stealth, thick, opacity=0.1] (a1) to[out=0,in=180] (a2);
        \draw[-stealth, thick] (a1) to[out=60,in=120] (a2);
        \draw[-stealth, thick, opacity=0.1] (a1) to[out=-60,in=-120] (a2);
        \draw[-stealth, thick, opacity=0.1] (a1) to[out=60,in=140] (a3);
        \draw[-stealth, thick, opacity=0.1] (a1) to[out=-60,in=-140] (a4);
        
        \draw[-stealth, thick, opacity=0.1] (a2) to[out=60,in=200] (a3);
        \draw[-stealth, thick, opacity=0.1] (a2) to[out=-60,in=-200] (a4);
        \draw[-stealth, thick, dashed] (a2) to[out=0,in=180] (a5);
        
        \draw[-stealth, thick, opacity=0.1] (a3) to[out=30,in=90] (a5);
        \draw[-stealth, thick, opacity=0.1] (a3) to[out=-30,in=135] (a5);
        
        \draw[-stealth, thick, opacity=0.1] (a4) to[out=-30,in=-90] (a5);
        \draw[-stealth, thick, opacity=0.1] (a4) to[out=30,in=-135] (a5);
    \end{scope} 
    \begin{scope}[xshift=90, yshift=-80, scale=0.9]
        \vertex[fill](a1) at (1,0) {};
        \vertex[fill](a2) at (2,0) {};
        \vertex[fill, opacity=0.1](a3) at (3,0.8) {};
        \vertex[fill](a4) at (3,-0.8) {};
        \vertex[fill](a5) at (4,0) {};
        
        \draw[-stealth, thick, opacity=0.1] (a1) to[out=0,in=180] (a2);
        \draw[-stealth, thick, opacity=0.1] (a1) to[out=60,in=120] (a2);
        \draw[-stealth, thick, dashed] (a1) to[out=-60,in=-120] (a2);
        \draw[-stealth, thick, opacity=0.1] (a1) to[out=60,in=140] (a3);
        \draw[-stealth, thick, opacity=0.1] (a1) to[out=-60,in=-140] (a4);
        
        \draw[-stealth, thick, opacity=0.1] (a2) to[out=60,in=200] (a3);
        \draw[-stealth, thick] (a2) to[out=-60,in=-200] (a4);
        \draw[-stealth, thick, opacity=0.1] (a2) to[out=0,in=180] (a5);
        
        \draw[-stealth, thick, opacity=0.1] (a3) to[out=30,in=90] (a5);
        \draw[-stealth, thick, opacity=0.1] (a3) to[out=-30,in=135] (a5);
        
        \draw[-stealth, thick,opacity=0.1] (a4) to[out=-30,in=-90] (a5);
        \draw[-stealth, thick] (a4) to[out=30,in=-135] (a5);
    \end{scope} 
    \begin{scope}[xshift=180, yshift=-80, scale=0.9]
        \vertex[fill](a1) at (1,0) {};
        \vertex[fill, opacity=0.1](a2) at (2,0) {};
        \vertex[fill, opacity=0.1](a3) at (3,0.8) {};
        \vertex[fill](a4) at (3,-0.8) {};
        \vertex[fill](a5) at (4,0) {};
        
        \draw[-stealth, thick, opacity=0.1] (a1) to[out=0,in=180] (a2);
        \draw[-stealth, thick, opacity=0.1] (a1) to[out=60,in=120] (a2);
        \draw[-stealth, thick, opacity=0.1] (a1) to[out=-60,in=-120] (a2);
        \draw[-stealth, thick, opacity=0.1] (a1) to[out=60,in=140] (a3);
        \draw[-stealth, thick, dashed] (a1) to[out=-60,in=-140] (a4);
        
        \draw[-stealth, thick, opacity=0.1] (a2) to[out=60,in=200] (a3);
        \draw[-stealth, thick, opacity=0.1] (a2) to[out=-60,in=-200] (a4);
        \draw[-stealth, thick, opacity=0.1] (a2) to[out=0,in=180] (a5);
        
        \draw[-stealth, thick, opacity=0.1] (a3) to[out=30,in=90] (a5);
        \draw[-stealth, thick, opacity=0.1] (a3) to[out=-30,in=135] (a5);
        
        \draw[-stealth, thick] (a4) to[out=-30,in=-90] (a5);
        \draw[-stealth, thick, opacity=0.1] (a4) to[out=30,in=-135] (a5);
    \end{scope}

\end{tikzpicture}

%% file: routeframing.tikz
\begin{tikzpicture}
\begin{scope}[xshift=0, yshift=0, scale=0.9]
        \vertex[fill](a1) at (1,0) {};
        \vertex[fill](a2) at (2,0) {};
        \vertex[fill](a3) at (3,0.8) {};
        \vertex[fill](a4) at (3,-0.8) {};
        \vertex[fill](a5) at (4,0) {};
        
        \draw[-stealth, thick] (a1) to[out=0,in=180] node[yshift=1.375mm] {\tiny $2$}  (a2);
        \draw[-stealth, thick] (a1) to[out=60,in=120] node[yshift=1.375mm] {\tiny $3$}  (a2);
        \draw[-stealth, thick] (a1) to[out=-60,in=-120] node[yshift=1.375mm] {\tiny $4$}  (a2);
        \draw[-stealth, thick] (a1) to[out=60,in=140] node[yshift=1.75mm] {\tiny $1$}  (a3);
        \draw[-stealth, thick] (a1) to[out=-60,in=-140]node[yshift=1.5mm] {\tiny $5$}  (a4);
        
        \draw[-stealth, thick] (a2) to[out=60,in=200] node[yshift=1.75mm] {\tiny $2$}  (a3);
        \draw[-stealth, thick] (a2) to[out=-60,in=-200]node[yshift=1.75mm] {\tiny $4$}  (a4);
        \draw[-stealth, thick] (a2) to[out=0,in=180] node[yshift=1.5mm] {\tiny $3$}  (a5);
        
        \draw[-stealth, thick] (a3) to[out=30,in=90] node[yshift=1.875mm] {\tiny $2$}  (a5);
        \draw[-stealth, thick] (a3) to[out=-30,in=135] node[yshift=1.75mm] {\tiny $1$} (a5);
        
        \draw[-stealth, thick] (a4) to[out=-30,in=-90] node[yshift=1.75mm] {\tiny $5$}  (a5);
        \draw[-stealth, thick] (a4) to[out=30,in=-135] node[yshift=1.875mm] {\tiny $4$}  (a5);
    
    \end{scope}   

    \begin{scope}[xshift=-180, yshift=-60, scale=0.8]
        \vertex[fill](a1) at (1,0) {};
        \vertex[fill,opacity=0.1](a2) at (2,0) {};%
        \vertex[fill](a3) at (3,0.8) {};
        \vertex[fill,opacity=0.1](a4) at (3,-0.8) {};%
        \vertex[fill](a5) at (4,0) {};
        
        \draw[-stealth, thick,opacity=0.1] (a1) to[out=0,in=180] (a2);%
        \draw[-stealth, thick,opacity=0.1] (a1) to[out=60,in=120] (a2);%
        \draw[-stealth, thick,opacity=0.1] (a1) to[out=-60,in=-120] (a2);%
        \draw[-stealth, thick] (a1) to[out=60,in=140] (a3);
        \draw[-stealth, thick,opacity=0.1] (a1) to[out=-60,in=-140] (a4);%
        
        \draw[-stealth, thick,opacity=0.1] (a2) to[out=60,in=200] (a3);%
        \draw[-stealth, thick,opacity=0.1] (a2) to[out=-60,in=-200] (a4);%
        \draw[-stealth, thick,opacity=0.1] (a2) to[out=0,in=180] (a5);%
        
        \draw[-stealth, thick,opacity=0.1] (a3) to[out=30,in=90] (a5);%
        \draw[-stealth, thick] (a3) to[out=-30,in=135] (a5);
        
        \draw[-stealth, thick,opacity=0.1] (a4) to[out=-30,in=-90] (a5);%
        \draw[-stealth, thick,opacity=0.1] (a4) to[out=30,in=-135] (a5);%

        \node at (4.5,0) {$\leq$}; 
        \node at (0.25,1) {$(\mathcal{R}, \leq ):$};
        \node at (2.5,-1.5) {$R_1$};
    \end{scope}  
    \begin{scope}[xshift=-90, yshift=-60, scale=0.8]
        \vertex[fill](a1) at (1,0) {};
        \vertex[fill](a2) at (2,0) {};
        \vertex[fill](a3) at (3,0.8) {};
        \vertex[fill, opacity=0.1](a4) at (3,-0.8) {};
        \vertex[fill](a5) at (4,0) {};
        
        \draw[-stealth, thick] (a1) to[out=0,in=180] (a2);
        \draw[-stealth, thick, opacity=0.1] (a1) to[out=60,in=120] (a2);%
        \draw[-stealth, thick, opacity=0.1] (a1) to[out=-60,in=-120] (a2);%
        \draw[-stealth, thick, opacity=0.1] (a1) to[out=60,in=140] (a3);%
        \draw[-stealth, thick, opacity=0.1] (a1) to[out=-60,in=-140] (a4);%
        
        \draw[-stealth, thick] (a2) to[out=60,in=200] (a3);
        \draw[-stealth, thick, opacity=0.1] (a2) to[out=-60,in=-200] (a4);
        \draw[-stealth, thick, opacity=0.1] (a2) to[out=0,in=180] (a5);
        
        \draw[-stealth, thick] (a3) to[out=30,in=90] (a5);
        \draw[-stealth, thick,opacity=0.1] (a3) to[out=-30,in=135] (a5);
        
        \draw[-stealth, thick, opacity=0.1] (a4) to[out=-30,in=-90] (a5);
        \draw[-stealth, thick, opacity=0.1] (a4) to[out=30,in=-135] (a5);

        \node at (4.5,0) {$\leq$}; 
        \node at (2.5,-1.5) {$R_2$};
    \end{scope} 
    \begin{scope}[xshift=0, yshift=-60, scale=0.8]
        \vertex[fill](a1) at (1,0) {};
        \vertex[fill](a2) at (2,0) {};
        \vertex[fill, opacity=0.1](a3) at (3,0.8) {};
        \vertex[fill, opacity=0.1](a4) at (3,-0.8) {};
        \vertex[fill](a5) at (4,0) {};
        
        \draw[-stealth, thick, opacity=0.1] (a1) to[out=0,in=180] (a2);
        \draw[-stealth, thick] (a1) to[out=60,in=120] (a2);
        \draw[-stealth, thick, opacity=0.1] (a1) to[out=-60,in=-120] (a2);
        \draw[-stealth, thick, opacity=0.1] (a1) to[out=60,in=140] (a3);
        \draw[-stealth, thick, opacity=0.1] (a1) to[out=-60,in=-140] (a4);
        
        \draw[-stealth, thick, opacity=0.1] (a2) to[out=60,in=200] (a3);
        \draw[-stealth, thick, opacity=0.1] (a2) to[out=-60,in=-200] (a4);
        \draw[-stealth, thick] (a2) to[out=0,in=180] (a5);
        
        \draw[-stealth, thick, opacity=0.1] (a3) to[out=30,in=90] (a5);
        \draw[-stealth, thick, opacity=0.1] (a3) to[out=-30,in=135] (a5);
        
        \draw[-stealth, thick, opacity=0.1] (a4) to[out=-30,in=-90] (a5);
        \draw[-stealth, thick, opacity=0.1] (a4) to[out=30,in=-135] (a5);

        \node at (4.5,0) {$\leq$}; 
        \node at (2.5,-1.5) {$R_3$};
    \end{scope} 
    \begin{scope}[xshift=90, yshift=-60, scale=0.8]
        \vertex[fill](a1) at (1,0) {};
        \vertex[fill](a2) at (2,0) {};
        \vertex[fill, opacity=0.1](a3) at (3,0.8) {};
        \vertex[fill](a4) at (3,-0.8) {};
        \vertex[fill](a5) at (4,0) {};
        
        \draw[-stealth, thick, opacity=0.1] (a1) to[out=0,in=180] (a2);
        \draw[-stealth, thick, opacity=0.1] (a1) to[out=60,in=120] (a2);
        \draw[-stealth, thick] (a1) to[out=-60,in=-120] (a2);
        \draw[-stealth, thick, opacity=0.1] (a1) to[out=60,in=140] (a3);
        \draw[-stealth, thick, opacity=0.1] (a1) to[out=-60,in=-140] (a4);
        
        \draw[-stealth, thick, opacity=0.1] (a2) to[out=60,in=200] (a3);
        \draw[-stealth, thick] (a2) to[out=-60,in=-200] (a4);
        \draw[-stealth, thick, opacity=0.1] (a2) to[out=0,in=180] (a5);
        
        \draw[-stealth, thick, opacity=0.1] (a3) to[out=30,in=90] (a5);
        \draw[-stealth, thick, opacity=0.1] (a3) to[out=-30,in=135] (a5);
        
        \draw[-stealth, thick,opacity=0.1] (a4) to[out=-30,in=-90] (a5);
        \draw[-stealth, thick] (a4) to[out=30,in=-135] (a5);

        \node at (4.5,0) {$\leq$}; 
        \node at (2.5,-1.5) {$R_4$};
    \end{scope} 
    \begin{scope}[xshift=180, yshift=-60, scale=0.8]
        \vertex[fill](a1) at (1,0) {};
        \vertex[fill, opacity=0.1](a2) at (2,0) {};
        \vertex[fill, opacity=0.1](a3) at (3,0.8) {};
        \vertex[fill](a4) at (3,-0.8) {};
        \vertex[fill](a5) at (4,0) {};
        
        \draw[-stealth, thick, opacity=0.1] (a1) to[out=0,in=180] (a2);
        \draw[-stealth, thick, opacity=0.1] (a1) to[out=60,in=120] (a2);
        \draw[-stealth, thick, opacity=0.1] (a1) to[out=-60,in=-120] (a2);
        \draw[-stealth, thick, opacity=0.1] (a1) to[out=60,in=140] (a3);
        \draw[-stealth, thick] (a1) to[out=-60,in=-140] (a4);
        
        \draw[-stealth, thick, opacity=0.1] (a2) to[out=60,in=200] (a3);
        \draw[-stealth, thick, opacity=0.1] (a2) to[out=-60,in=-200] (a4);
        \draw[-stealth, thick, opacity=0.1] (a2) to[out=0,in=180] (a5);
        
        \draw[-stealth, thick, opacity=0.1] (a3) to[out=30,in=90] (a5);
        \draw[-stealth, thick, opacity=0.1] (a3) to[out=-30,in=135] (a5);
        
        \draw[-stealth, thick] (a4) to[out=-30,in=-90] (a5);
        \draw[-stealth, thick, opacity=0.1] (a4) to[out=30,in=-135] (a5);

        \node at (2.5,-1.5) {$R_5$};
    \end{scope} 
\end{tikzpicture}

%% file: notDKK.tikz
\begin{tikzpicture}
        \vertex[fill](a1) at (1,0) {};
        \vertex[fill](a2) at (2,0) {};
        \vertex[fill](a3) at (3,0) {};
        \vertex[fill](a4) at (4,0) {};

        
        \draw[-stealth, thick] (a1) to[out=60,in=120] node[yshift=0.75mm, xshift=5mm] {\tiny $1$}  (a3);
        \draw[-stealth, thick] (a1) to[out=60,in=120] node[yshift=1mm,xshift=3mm] {\tiny $2$} (a2);
        \draw[-stealth, thick] (a1) to[out=0,in=180] node[yshift=1.5mm] {\tiny $3$}  (a2);

        \draw[-stealth, thick] (a2) to[out=0,in=180] node[yshift=1.5mm] {\tiny $2$}  (a3);
        \draw[-stealth, thick] (a2) to[out=-60,in=-120] node[yshift=1.5mm] {\tiny $3$}  (a4);

        \draw[-stealth, thick] (a3) to[out=0,in=180] node[yshift=1.5mm] {\tiny $1$} (a4);
        \draw[-stealth, thick] (a3) to[out=-60,in=-120] node[yshift=1.5mm] {\tiny $2$}  (a4);
\end{tikzpicture}

%% file: quotientexample.tikz
\begin{tikzpicture}
        \vertex[fill](a1) at (1,0) {};
        \vertex[fill](a2) at (2,0) {};
        \vertex[fill](a3) at (3,0) {};
        \vertex[fill](a4) at (4,0) {};

        \node[below, color=orange] at (a1) {\small $s$};
        \node[below, color=orange] at (2,0) {\small $1$};
        \node[above, color=orange] at (3,0) {\small $2$};
        \node[above, color=orange] at (4,0) {\small $t$};
        
        \draw[-stealth, thick] (a1) to[out=60,in=120] node[yshift=0.75mm, xshift=5mm] {\tiny $1$}  (a3);
        \draw[-stealth, thick] (a1) to[out=60,in=120] node[yshift=1mm,xshift=3mm] {\tiny $2$} (a2);
        \draw[-stealth, thick] (a1) to[out=0,in=180] node[yshift=1.5mm] {\tiny $3$}  (a2);

        \draw[-stealth, thick] (a2) to[out=0,in=180] node[yshift=1.5mm] {\tiny $2$}  (a3);
        \draw[-stealth, thick] (a2) to[out=-60,in=-120] node[yshift=1.5mm] {\tiny $3$}  (a4);

        \draw[-stealth, thick] (a3) to[out=0,in=180] node[yshift=1.5mm] {\tiny $1$} (a4);
        \draw[-stealth, thick] (a3) to[out=-60,in=-120] node[yshift=1.5mm] {\tiny $2$}  (a4);
\end{tikzpicture}

\begin{tikzpicture}%
	[scale=1.000000,
	back/.style={loosely dotted, thick},
	edge/.style={color=green},
	facet/.style={fill=green,fill opacity=0.3}]
    
    \coordinate (0,0) at (0,0);
    \coordinate (0,2) at (0,2);
    \coordinate (2,0) at (2,0);
    \coordinate (2,2) at (2,2);
    \coordinate (3,1) at (3,1);
    \fill[facet] (3, 1) -- (2, 0) -- (0, 0) -- (0, 2) -- (2, 2) -- cycle {};
    \draw[edge] (0, 0) -- (0, 2);
    \draw[edge] (0, 0) -- (2, 0);
    \draw[edge] (0, 2) -- (2, 2);
    \draw[edge] (2, 0) -- (3, 1);
    \draw[edge] (2, 2) -- (3, 1);
    \node at (1,1) [circle,fill,inner sep=0.5pt, color=green]{};
    \node[xshift=-15, yshift=-17] at (0, 0)     {
    \begin{tikzpicture}[scale=0.6666]
            \vertex[fill](a1) at (1,0) {};
            \vertex[fill](a2) at (2,0) {};
            \vertex[fill,opacity=0.1](a3) at (3,0) {};
            \vertex[fill](a4) at (4,0) {};
            
            \draw[-stealth, thick, opacity=0.1] (a1) to[out=60,in=120] node[yshift=0.75mm, xshift=5mm] {\tiny $1$}  (a3);
            \draw[-stealth, thick] (a1) to[out=60,in=120] node[yshift=1mm,xshift=3mm] {\tiny $2$} (a2);
            \draw[-stealth, thick, opacity=0.1] (a1) to[out=0,in=180] node[yshift=1.5mm] {\tiny $3$}  (a2);
    
            \draw[-stealth, thick,opacity=0.1] (a2) to[out=0,in=180] node[yshift=1.5mm] {\tiny $2$}  (a3);
            \draw[-stealth, thick] (a2) to[out=-60,in=-120] node[yshift=1.5mm] {\tiny $3$}  (a4);
    
            \draw[-stealth, thick, opacity=0.1] (a3) to[out=0,in=180] node[yshift=1.5mm] {\tiny $1$} (a4);
            \draw[-stealth, thick, opacity=0.1] (a3) to[out=-60,in=-120] node[yshift=1.5mm] {\tiny $2$}  (a4);

            \node[] at (0.5,0) {$\varphi \bigl($};
            \node[] at (4.25,0) {$\bigr)$};
    \end{tikzpicture}
    };
    \node[xshift=-15, yshift=17] at (0, 2)     {
    \begin{tikzpicture}[scale=0.6666]
            \vertex[fill](a1) at (1,0) {};
            \vertex[fill, opacity=0.1](a2) at (2,0) {};
            \vertex[fill](a3) at (3,0) {};
            \vertex[fill](a4) at (4,0) {};
            
            \draw[-stealth, thick] (a1) to[out=60,in=120] node[yshift=0.75mm, xshift=5mm] {\tiny $1$}  (a3);
            \draw[-stealth, thick, opacity=0.1] (a1) to[out=60,in=120] node[yshift=1mm,xshift=3mm] {\tiny $2$} (a2);
            \draw[-stealth, thick, opacity=0.1] (a1) to[out=0,in=180] node[yshift=1.5mm] {\tiny $3$}  (a2);
    
            \draw[-stealth, thick, opacity=0.1] (a2) to[out=0,in=180] node[yshift=1.5mm] {\tiny $2$}  (a3);
            \draw[-stealth, thick, opacity=0.1] (a2) to[out=-60,in=-120] node[yshift=1.5mm] {\tiny $3$}  (a4);
    
            \draw[-stealth, thick, opacity=0.1] (a3) to[out=0,in=180] node[yshift=1.5mm] {\tiny $1$} (a4);
            \draw[-stealth, thick] (a3) to[out=-60,in=-120] node[yshift=1.5mm] {\tiny $2$}  (a4);

            \node[] at (0.5,0) {$\varphi \bigl($};
            \node[] at (4.25,0) {$\bigr)$};
    \end{tikzpicture}
    };
    \node[xshift=15,yshift=-17] at (2, 0)     {
    \begin{tikzpicture}[scale=0.6666]
            \vertex[fill](a1) at (1,0) {};
            \vertex[fill](a2) at (2,0) {};
            \vertex[fill](a3) at (3,0) {};
            \vertex[fill](a4) at (4,0) {};
            
            \draw[-stealth, thick, opacity=0.1] (a1) to[out=60,in=120] node[yshift=0.75mm, xshift=5mm] {\tiny $1$}  (a3);
            \draw[-stealth, thick ] (a1) to[out=60,in=120] node[yshift=1mm,xshift=3mm] {\tiny $2$} (a2);
            \draw[-stealth, thick, opacity=0.1] (a1) to[out=0,in=180] node[yshift=1.5mm] {\tiny $3$}  (a2);
    
            \draw[-stealth, thick] (a2) to[out=0,in=180] node[yshift=1.5mm] {\tiny $2$}  (a3);
            \draw[-stealth, thick, opacity=0.1] (a2) to[out=-60,in=-120] node[yshift=1.5mm] {\tiny $3$}  (a4);
    
            \draw[-stealth, thick] (a3) to[out=0,in=180] node[yshift=1.5mm] {\tiny $1$} (a4);
            \draw[-stealth, thick, opacity=0.1] (a3) to[out=-60,in=-120] node[yshift=1.5mm] {\tiny $2$}  (a4);

            \node[] at (0.5,0) {$\varphi \bigl($};
            \node[] at (4.25,0) {$\bigr)$};
    \end{tikzpicture}
    };
    \node[xshift=15,yshift=17] at (2, 2)     {
    \begin{tikzpicture}[scale=0.6666]
            \vertex[fill](a1) at (1,0) {};
            \vertex[fill](a2) at (2,0) {};
            \vertex[fill](a3) at (3,0) {};
            \vertex[fill](a4) at (4,0) {};
            
            \draw[-stealth, thick, opacity=0.1] (a1) to[out=60,in=120] node[yshift=0.75mm, xshift=5mm] {\tiny $1$}  (a3);
            \draw[-stealth, thick, opacity=0.1 ] (a1) to[out=60,in=120] node[yshift=1mm,xshift=3mm] {\tiny $2$} (a2);
            \draw[-stealth, thick] (a1) to[out=0,in=180] node[yshift=1.5mm] {\tiny $3$}  (a2);
    
            \draw[-stealth, thick] (a2) to[out=0,in=180] node[yshift=1.5mm] {\tiny $2$}  (a3);
            \draw[-stealth, thick, opacity=0.1] (a2) to[out=-60,in=-120] node[yshift=1.5mm] {\tiny $3$}  (a4);
    
            \draw[-stealth, thick, opacity=0.1] (a3) to[out=0,in=180] node[yshift=1.5mm] {\tiny $1$} (a4);
            \draw[-stealth, thick] (a3) to[out=-60,in=-120] node[yshift=1.5mm] {\tiny $2$}  (a4);

            \node[] at (0.5,0) {$\varphi \bigl($};
            \node[] at (4.25,0) {$\bigr)$};
    \end{tikzpicture}
    };
    \node[xshift=45, yshift=3] at (3, 1)     {
    \begin{tikzpicture}[scale=0.6666]
            \vertex[fill](a1) at (1,0) {};
            \vertex[fill](a2) at (2,0) {};
            \vertex[fill](a3) at (3,0) {};
            \vertex[fill](a4) at (4,0) {};
            
            \draw[-stealth, thick, opacity=0.1] (a1) to[out=60,in=120] node[yshift=0.75mm, xshift=5mm] {\tiny $1$}  (a3);
            \draw[-stealth, thick, opacity=0.1 ] (a1) to[out=60,in=120] node[yshift=1mm,xshift=3mm] {\tiny $2$} (a2);
            \draw[-stealth, thick] (a1) to[out=0,in=180] node[yshift=1.5mm] {\tiny $3$}  (a2);
    
            \draw[-stealth, thick] (a2) to[out=0,in=180] node[yshift=1.5mm] {\tiny $2$}  (a3);
            \draw[-stealth, thick, opacity=0.1] (a2) to[out=-60,in=-120] node[yshift=1.5mm] {\tiny $3$}  (a4);
    
            \draw[-stealth, thick] (a3) to[out=0,in=180] node[yshift=1.5mm] {\tiny $1$} (a4);
            \draw[-stealth, thick, opacity=0.1] (a3) to[out=-60,in=-120] node[yshift=1.5mm] {\tiny $2$}  (a4);

            \node[] at (0.5,0) {$\varphi \bigl($};
            \node[] at (4.25,0) {$\bigr)$};
    \end{tikzpicture}
    };
\end{tikzpicture}
\begin{tikzpicture}%
	[scale=1.000000,
	back/.style={loosely dotted, thick},
	edge/.style={color=green},
	facet/.style={fill=green,fill opacity=0.3}]
    
    \coordinate (0,0) at (0,0);
    \coordinate (0,2) at (0,2);
    \coordinate (2,0) at (2,0);
    \coordinate (2,2) at (2,2);
    \coordinate (3,1) at (3,1);
    \fill[facet] (3, 1) -- (2, 0) -- (0, 0) -- (0, 2) -- (2, 2) -- cycle {};
    \draw[edge] (0, 0) -- (0, 2);
    \draw[edge] (0, 0) -- (2, 0);
    \draw[edge] (0, 2) -- (2, 2);
    \draw[edge] (2, 0) -- (3, 1);
    \draw[edge] (2, 2) -- (3, 1);
    \node at (1,1) [circle,fill,inner sep=0.5pt, color=green]{};
    \node[xshift=-15, yshift=-17] at (0, 0)     {$e^1_2 - e^1_3$};
    \node[xshift=-15, yshift=17] at (0, 2)     {$e^2_1 - e^2_2$};
    \node[xshift=15,yshift=-17] at (2, 0)     {$e^2_2 - e^2_1$};
    \node[xshift=15,yshift=17] at (2, 2)     {$e^1_3-e^1_2$};
    \node[xshift=45, yshift=3] at (3, 1)     {$e^1_3-e^1_2+e_2^2-e^2_1$};
\end{tikzpicture}

%% file: quotient3dexample.tikz
\begin{tikzpicture}
        \vertex[fill](a1) at (1,0) {};
        \vertex[fill](a2) at (2.5,0) {};
        \vertex[fill](a3) at (4,0) {};
        
        \draw[-stealth, thick] (a1) to[out=80,in=100] node[yshift=1.5mm] {\tiny $4$} (a2);
        \draw[-stealth, thick] (a1) to[out=30,in=150] node[yshift=1.25mm] {\tiny $3$} (a2);
        \draw[-stealth, thick] (a1) to[out=-30,in=-150] node[yshift=1.25mm] {\tiny $2$} (a2);
        \draw[-stealth, thick] (a1) to[out=-80,in=-100] node[yshift=1.25mm] {\tiny $1$} (a2);

        \draw[-stealth, thick] (a2) to[out=80,in=100] node[yshift=1.5mm] {\tiny $4$} (a3);
        \draw[-stealth, thick] (a2) to[out=30,in=150] node[yshift=1.25mm] {\tiny $3$} (a3);
        \draw[-stealth, thick] (a2) to[out=-30,in=-150] node[yshift=1.25mm] {\tiny $2$} (a3);
        \draw[-stealth, thick] (a2) to[out=-80,in=-100] node[yshift=1.25mm] {\tiny $1$} (a3);   
\end{tikzpicture}
\:
\begin{tikzpicture}%
	[x={(-0.022585cm, -0.481709cm)},
	y={(0.999745cm, -0.010857cm)},
	z={(-0.000029cm, 0.876264cm)},
	scale=1.000000,
    vertex/.style={opacity=1},
	back/.style={loosely dotted, thick},
	edge/.style={color=green},
	facet/.style={fill=green,fill opacity=0.300000}]
    %
    %
    
    \coordinate (0.70711, -0.40825, 2.30940) at (0.70711, -0.40825, 2.30940);
    \coordinate (1.41421, 0.81650, 2.30940) at (1.41421, 0.81650, 2.30940);
    \coordinate (0.00000, 0.81650, 2.30940) at (0.00000, 0.81650, 2.30940);
    \coordinate (1.41421, -0.81650, 1.15470) at (1.41421, -0.81650, 1.15470);
    \coordinate (2.12132, 0.40825, 1.15470) at (2.12132, 0.40825, 1.15470);
    \coordinate (0.00000, -0.81650, 1.15470) at (0.00000, -0.81650, 1.15470);
    \coordinate (0.70711, 1.22474, 0.00000) at (0.70711, 1.22474, 0.00000);
    \coordinate (1.41421, 1.63299, 1.15470) at (1.41421, 1.63299, 1.15470);
    \coordinate (-0.70711, 0.40825, 1.15470) at (-0.70711, 0.40825, 1.15470);
    \coordinate (0.00000, 1.63299, 1.15470) at (0.00000, 1.63299, 1.15470);
    \coordinate (1.41421, 0.00000, 0.00000) at (1.41421, 0.00000, 0.00000);
    \coordinate (0.00000, 0.00000, 0.00000) at (0.00000, 0.00000, 0.00000);
    \draw[edge,back] (0.00000, 0.81650, 2.30940) -- (-0.70711, 0.40825, 1.15470);
    \draw[edge,back] (0.00000, -0.81650, 1.15470) -- (-0.70711, 0.40825, 1.15470);
    \draw[edge,back] (0.00000, -0.81650, 1.15470) -- (0.00000, 0.00000, 0.00000);
    \draw[edge,back] (0.70711, 1.22474, 0.00000) -- (0.00000, 1.63299, 1.15470);
    \draw[edge,back] (0.70711, 1.22474, 0.00000) -- (0.00000, 0.00000, 0.00000);
    \draw[edge,back] (-0.70711, 0.40825, 1.15470) -- (0.00000, 1.63299, 1.15470);
    \draw[edge,back] (-0.70711, 0.40825, 1.15470) -- (0.00000, 0.00000, 0.00000);
    \draw[edge,back] (1.41421, 0.00000, 0.00000) -- (0.00000, 0.00000, 0.00000);
    \node[vertex] at (-0.70711, 0.40825, 1.15470)     {};
    \node[vertex] at (0.00000, 0.00000, 0.00000)     {};
    \fill[facet] (0.70711, 1.22474, 0.00000) -- (1.41421, 0.00000, 0.00000) -- (2.12132, 0.40825, 1.15470) -- (1.41421, 1.63299, 1.15470) -- cycle {};
    \fill[facet] (2.12132, 0.40825, 1.15470) -- (1.41421, 0.81650, 2.30940) -- (0.70711, -0.40825, 2.30940) -- (1.41421, -0.81650, 1.15470) -- cycle {};
    \fill[facet] (0.00000, -0.81650, 1.15470) -- (0.70711, -0.40825, 2.30940) -- (1.41421, -0.81650, 1.15470) -- cycle {};
    \fill[facet] (1.41421, 1.63299, 1.15470) -- (1.41421, 0.81650, 2.30940) -- (2.12132, 0.40825, 1.15470) -- cycle {};
    \fill[facet] (0.00000, 0.81650, 2.30940) -- (0.70711, -0.40825, 2.30940) -- (1.41421, 0.81650, 2.30940) -- cycle {};
    \fill[facet] (0.00000, 1.63299, 1.15470) -- (0.00000, 0.81650, 2.30940) -- (1.41421, 0.81650, 2.30940) -- (1.41421, 1.63299, 1.15470) -- cycle {};
    \fill[facet] (1.41421, 0.00000, 0.00000) -- (1.41421, -0.81650, 1.15470) -- (2.12132, 0.40825, 1.15470) -- cycle {};
    \draw[edge] (0.70711, -0.40825, 2.30940) -- (1.41421, 0.81650, 2.30940);
    \draw[edge] (0.70711, -0.40825, 2.30940) -- (0.00000, 0.81650, 2.30940);
    \draw[edge] (0.70711, -0.40825, 2.30940) -- (1.41421, -0.81650, 1.15470);
    \draw[edge] (0.70711, -0.40825, 2.30940) -- (0.00000, -0.81650, 1.15470);
    \draw[edge] (1.41421, 0.81650, 2.30940) -- (0.00000, 0.81650, 2.30940);
    \draw[edge] (1.41421, 0.81650, 2.30940) -- (2.12132, 0.40825, 1.15470);
    \draw[edge] (1.41421, 0.81650, 2.30940) -- (1.41421, 1.63299, 1.15470);
    \draw[edge] (0.00000, 0.81650, 2.30940) -- (0.00000, 1.63299, 1.15470);
    \draw[edge] (1.41421, -0.81650, 1.15470) -- (2.12132, 0.40825, 1.15470);
    \draw[edge] (1.41421, -0.81650, 1.15470) -- (0.00000, -0.81650, 1.15470);
    \draw[edge] (1.41421, -0.81650, 1.15470) -- (1.41421, 0.00000, 0.00000);
    \draw[edge] (2.12132, 0.40825, 1.15470) -- (1.41421, 1.63299, 1.15470);
    \draw[edge] (2.12132, 0.40825, 1.15470) -- (1.41421, 0.00000, 0.00000);
    \draw[edge] (0.70711, 1.22474, 0.00000) -- (1.41421, 1.63299, 1.15470);
    \draw[edge] (0.70711, 1.22474, 0.00000) -- (1.41421, 0.00000, 0.00000);
    \draw[edge] (1.41421, 1.63299, 1.15470) -- (0.00000, 1.63299, 1.15470);
    \node[vertex] at (0.70711, -0.40825, 2.30940)     {};
    \node[vertex] at (1.41421, 0.81650, 2.30940)     {};
    \node[vertex] at (0.00000, 0.81650, 2.30940)     {};
    \node[vertex] at (1.41421, -0.81650, 1.15470)     {};
    \node[vertex] at (2.12132, 0.40825, 1.15470)     {};
    \node[vertex] at (0.00000, -0.81650, 1.15470)     {};
    \node[vertex] at (0.70711, 1.22474, 0.00000)     {};
    \node[vertex] at (1.41421, 1.63299, 1.15470)     {};
    \node[vertex] at (0.00000, 1.63299, 1.15470)     {};
    \node[vertex] at (1.41421, 0.00000, 0.00000)     {};
    \node at (0.707106,0.40825,1.15470) [circle,fill,inner sep=0.5pt, color=green]{};
\end{tikzpicture}
\hspace{0.25in}
\begin{tikzpicture}
        \vertex[fill](a1) at (1,0) {};
        \vertex[fill](a2) at (2,0) {};
        \vertex[fill](a3) at (3,0) {};
        \vertex[fill](a4) at (4,0) {};
        \vertex[fill](a5) at (5,0) {};
        
        \draw[-stealth, thick] (a1) to[out=60,in=120] node[yshift=1.5mm] {\tiny $1$}  (a4); 
        \draw[-stealth, thick] (a1) to[out=60,in=120] node[yshift=0.75mm, xshift=5mm] {\tiny $2$}  (a3);
        \draw[-stealth, thick] (a1) to[out=60,in=120] node[yshift=1mm,xshift=3mm] {\tiny $3$} (a2);
        \draw[-stealth, thick] (a1) to[out=0,in=180] node[yshift=1.5mm] {\tiny $4$}  (a2);

        \draw[-stealth, thick] (a2) to[out=0,in=180] node[yshift=1.5mm] {\tiny $3$}  (a3);
        \draw[-stealth, thick] (a2) to[out=-60,in=-120] node[yshift=1.5mm] {\tiny $4$}  (a5);

        \draw[-stealth, thick] (a3) to[out=0,in=180] node[yshift=1.5mm] {\tiny $2$} (a4);
        \draw[-stealth, thick] (a3) to[out=-60,in=-120] node[yshift=1.5mm] {\tiny $3$}  (a5);

        \draw[-stealth, thick] (a4) to[out=0,in=180] node[yshift=1.5mm] {\tiny $1$}  (a5);
        \draw[-stealth, thick] (a4) to[out=-60,in=-120] node[yshift=1.25mm] {\tiny $2$} (a5);

\end{tikzpicture}
\:
\begin{tikzpicture}%
	[x={(-0.022585cm, -0.481709cm)},
	y={(0.999745cm, -0.010857cm)},
	z={(-0.000029cm, 0.876264cm)},
	scale=1.000000,
    vertex/.style={opacity=1},
	back/.style={loosely dotted, thick},
	edge/.style={color=green},
	facet/.style={fill=green,fill opacity=0.300000}]
    %
    %
    
    \coordinate (1.00000, 0.57735, 3.26599) at (1.00000, 0.57735, 3.26599);
    \coordinate (1.00000, 1.73205, 2.44949) at (1.00000, 1.73205, 2.44949);
    \coordinate (2.00000, 1.15470, 1.63299) at (2.00000, 1.15470, 1.63299);
    \coordinate (0.00000, 0.00000, 2.44949) at (0.00000, 0.00000, 2.44949);
    \coordinate (0.00000, 1.15470, 1.63299) at (0.00000, 1.15470, 1.63299);
    \coordinate (1.00000, -0.57735, 1.63299) at (1.00000, -0.57735, 1.63299);
    \coordinate (1.00000, 1.73205, 0.00000) at (1.00000, 1.73205, 0.00000);
    \coordinate (2.00000, 0.00000, 0.00000) at (2.00000, 0.00000, 0.00000);
    \coordinate (0.00000, 0.00000, 0.00000) at (0.00000, 0.00000, 0.00000);
    \draw[edge,back] (1.00000, 1.73205, 2.44949) -- (0.00000, 1.15470, 1.63299);
    \draw[edge,back] (0.00000, 0.00000, 2.44949) -- (0.00000, 1.15470, 1.63299);
    \draw[edge,back] (0.00000, 0.00000, 2.44949) -- (0.00000, 0.00000, 0.00000);
    \draw[edge,back] (0.00000, 1.15470, 1.63299) -- (1.00000, 1.73205, 0.00000);
    \draw[edge,back] (0.00000, 1.15470, 1.63299) -- (0.00000, 0.00000, 0.00000);
    \draw[edge,back] (1.00000, -0.57735, 1.63299) -- (0.00000, 0.00000, 0.00000);
    \draw[edge,back] (1.00000, 1.73205, 0.00000) -- (0.00000, 0.00000, 0.00000);
    \draw[edge,back] (2.00000, 0.00000, 0.00000) -- (0.00000, 0.00000, 0.00000);
    \node[vertex] at (0.00000, 1.15470, 1.63299)     {};
    \node[vertex] at (0.00000, 0.00000, 0.00000)     {};
    \fill[facet] (2.00000, 0.00000, 0.00000) -- (2.00000, 1.15470, 1.63299) -- (1.00000, 1.73205, 0.00000) -- cycle {};
    \fill[facet] (2.00000, 0.00000, 0.00000) -- (2.00000, 1.15470, 1.63299) -- (1.00000, 0.57735, 3.26599) -- (1.00000, -0.57735, 1.63299) -- cycle {};
    \fill[facet] (1.00000, -0.57735, 1.63299) -- (1.00000, 0.57735, 3.26599) -- (0.00000, 0.00000, 2.44949) -- cycle {};
    \fill[facet] (2.00000, 1.15470, 1.63299) -- (1.00000, 0.57735, 3.26599) -- (1.00000, 1.73205, 2.44949) -- cycle {};
    \fill[facet] (1.00000, 1.73205, 0.00000) -- (1.00000, 1.73205, 2.44949) -- (2.00000, 1.15470, 1.63299) -- cycle {};
    \draw[edge] (1.00000, 0.57735, 3.26599) -- (1.00000, 1.73205, 2.44949);
    \draw[edge] (1.00000, 0.57735, 3.26599) -- (2.00000, 1.15470, 1.63299);
    \draw[edge] (1.00000, 0.57735, 3.26599) -- (0.00000, 0.00000, 2.44949);
    \draw[edge] (1.00000, 0.57735, 3.26599) -- (1.00000, -0.57735, 1.63299);
    \draw[edge] (1.00000, 1.73205, 2.44949) -- (2.00000, 1.15470, 1.63299);
    \draw[edge] (1.00000, 1.73205, 2.44949) -- (1.00000, 1.73205, 0.00000);
    \draw[edge] (2.00000, 1.15470, 1.63299) -- (1.00000, 1.73205, 0.00000);
    \draw[edge] (2.00000, 1.15470, 1.63299) -- (2.00000, 0.00000, 0.00000);
    \draw[edge] (0.00000, 0.00000, 2.44949) -- (1.00000, -0.57735, 1.63299);
    \draw[edge] (1.00000, -0.57735, 1.63299) -- (2.00000, 0.00000, 0.00000);
    \draw[edge] (1.00000, 1.73205, 0.00000) -- (2.00000, 0.00000, 0.00000);
    \node[vertex] at (1.00000, 0.57735, 3.26599)     {};
    \node[vertex] at (1.00000, 1.73205, 2.44949)     {};
    \node[vertex] at (2.00000, 1.15470, 1.63299)     {};
    \node[vertex] at (0.00000, 0.00000, 2.44949)     {};
    \node[vertex] at (1.00000, -0.57735, 1.63299)     {};
    \node[vertex] at (1.00000, 1.73205, 0.00000)     {};
    \node[vertex] at (2.00000, 0.00000, 0.00000)     {};

\end{tikzpicture}

%% file: routeequation.tikz
\begin{tikzpicture}
    \begin{scope}[xshift=-120, yshift=0, scale=0.9]
        \vertex[fill](a1) at (1,0) {};
        \vertex[fill](a2) at (2,0) {};
        \vertex[fill](a3) at (3,0.8) {};
        \vertex[fill,](a4) at (3,-0.8) {};
        \vertex[fill](a5) at (4,0) {};

        \node[left, color=orange] at (a1) {\small $s$};
        \node[below, color=orange] at (a2) {\small $1$};
        \node[below, color=orange] at (a3) {\small $2$};
        \node[below, color=orange] at (a4) {\small $3$};
        \node[right, color=orange] at (a5) {\small $t$};
        
        \draw[-stealth, thick] (a1) to[out=0,in=180] node[yshift=1.375mm] {\tiny $2$}  (a2);
        \draw[-stealth, thick] (a1) to[out=60,in=120] node[yshift=1.375mm] {\tiny $3$}  (a2);
        \draw[-stealth, thick,dashed] (a1) to[out=-60,in=-120] node[yshift=1.375mm] {\tiny $4$}  (a2);
        \draw[-stealth, thick] (a1) to[out=60,in=140] node[yshift=1.75mm] {\tiny $1$}  (a3);
        \draw[-stealth, thick,dashed] (a1) to[out=-60,in=-140]node[yshift=1.5mm] {\tiny $5$}  (a4);
        
        \draw[-stealth, thick,dashed] (a2) to[out=60,in=200] node[yshift=1.75mm] {\tiny $2$}  (a3);
        \draw[-stealth, thick] (a2) to[out=-60,in=-200]node[yshift=1.75mm] {\tiny $4$}  (a4);
        \draw[-stealth, thick,dashed] (a2) to[out=0,in=180] node[yshift=1.5mm] {\tiny $3$}  (a5);
        
        \draw[-stealth, thick] (a3) to[out=30,in=90] node[yshift=1.875mm] {\tiny $2$}  (a5);
        \draw[-stealth, thick,dashed] (a3) to[out=-30,in=135] node[yshift=1.75mm] {\tiny $1$} (a5);
        
        \draw[-stealth, thick] (a4) to[out=-30,in=-90] node[yshift=1.75mm] {\tiny $5$}  (a5);
        \draw[-stealth, thick] (a4) to[out=30,in=-135] node[yshift=1.875mm] {\tiny $4$}  (a5);
    
    \end{scope} 
    \begin{scope}[xshift=0, yshift=0, scale=0.9]
        \vertex[fill](a1) at (1,0) {};
        \vertex[fill](a2) at (2,0) {};
        \vertex[fill](a3) at (3,0.8) {};
        \vertex[fill,opacity=0.1](a4) at (3,-0.8) {};
        \vertex[fill](a5) at (4,0) {};

        \node[left, color=orange] at (a1) {\small $s$};
        \node[below, color=orange] at (a2) {\small $1$};
        \node[below, color=orange] at (a3) {\small $2$};
        \node[below, color=orange, opacity=0.1] at (a4) {$3$};
        \node[right, color=orange] at (a5) {\small $t$};
        
        \draw[-stealth, thick,opacity=0.1] (a1) to[out=0,in=180] node[yshift=1.375mm] {\tiny $2$}  (a2);
        \draw[-stealth, thick] (a1) to[out=60,in=120] node[yshift=1.375mm] {\tiny $3$}  (a2);
        \draw[-stealth, thick,dashed,opacity=0.1] (a1) to[out=-60,in=-120] node[yshift=1.375mm] {\tiny $4$}  (a2);
        \draw[-stealth, thick,opacity=0.1] (a1) to[out=60,in=140] node[yshift=1.75mm] {\tiny $1$}  (a3);
        \draw[-stealth, thick,dashed,opacity=0.1] (a1) to[out=-60,in=-140]node[yshift=1.5mm] {\tiny $5$}  (a4);
        
        \draw[-stealth, thick,dashed] (a2) to[out=60,in=200] node[yshift=1.75mm] {\tiny $2$}  (a3);
        \draw[-stealth, thick,opacity=0.1] (a2) to[out=-60,in=-200]node[yshift=1.75mm] {\tiny $4$}  (a4);
        \draw[-stealth, thick,dashed,opacity=0.1] (a2) to[out=0,in=180] node[yshift=1.5mm] {\tiny $3$}  (a5);
        
        \draw[-stealth, thick, opacity=0.1] (a3) to[out=30,in=90] node[yshift=1.875mm] {\tiny $2$}  (a5);
        \draw[-stealth, thick,dashed] (a3) to[out=-30,in=135] node[yshift=1.75mm] {\tiny $1$} (a5);
        
        \draw[-stealth, thick,opacity=0.1] (a4) to[out=-30,in=-90] node[yshift=1.75mm] {\tiny $5$}  (a5);
        \draw[-stealth, thick,opacity=0.1] (a4) to[out=30,in=-135] node[yshift=1.875mm] {\tiny $4$}  (a5);
    
    \end{scope}   

\end{tikzpicture}

%% file: stronglyplanar.tikz
\begin{tikzpicture}[scale=0.9]
    \draw[->, thick, gray] (2,1) to[out=60,in=120]  (9,1); 
    \draw[<-, thick, gray] (2,-0.5) to[out=-60,in=-120] (9,-0.5); 

    \begin{scope}[xshift=0, yshift=0, scale=1]
        
        \vertex[fill](a1) at (0.5,0) {};
        \vertex[fill](a2) at (1.5,0) {};
        \vertex[fill](a3) at (2,0.5) {};
        \vertex[fill](a4) at (3,0.5) {};

        \node[] at (0.5,0.75) {$G$};
        
        \draw[-stealth, thick] (a1) to[out=60,in=120] node[yshift=0.75mm, xshift=5mm] {}  (a3);
        \draw[-stealth, thick] (a1) to[out=60,in=120] node[yshift=1mm,xshift=3mm] {} (a2);
        \draw[-stealth, thick] (a1) to[out=-60,in=-120] node[yshift=1.5mm] {}  (a2);

        \draw[-stealth, thick] (a2) --  (a3);
        \draw[-stealth, thick] (a2) to[out=-60,in=-120] node[yshift=1.5mm] {}  (a4);

        \draw[-stealth, thick] (a3) to[out=60,in=120] node[yshift=1.5mm] {} (a4);
        \draw[-stealth, thick] (a3) to[out=-60,in=-120] node[yshift=1.5mm] {}  (a4);
    \end{scope}

    \begin{scope}[xshift=107, yshift=65, scale=1]
        \fill[white] (0,0) rectangle (3.5,1);
    
        \vertex[fill,cyan](a1) at (1,0) {};
        \vertex[fill,cyan](a2) at (1.5,0.5) {};{};
        \vertex[fill,cyan](a3) at (2,0) {};
        \vertex[fill,cyan](a4) at (2.5,0.5) {};

        \draw[thick,cyan] (a1) -- (a2);
        \draw[thick,cyan] (a3) -- (a2);
        \draw[thick,cyan] (a3) -- (a4);

    \end{scope}

   \begin{scope}[xshift=107, yshift=65, scale=1]
        \vertex[fill](a1) at (0.5,0) {};
        \vertex[fill](a2) at (1.5,0) {};
        \vertex[fill](a3) at (2,0.5) {};
        \vertex[fill](a4) at (3,0.5) {};

        
        \draw[-stealth, thick] (a1) to[out=60,in=120] node[yshift=0.75mm, xshift=5mm] {}  (a3);
        \draw[-stealth, thick] (a1) to[out=60,in=120] node[yshift=1mm,xshift=3mm] {} (a2);
        \draw[-stealth, thick] (a1) to[out=-60,in=-120] node[yshift=1.5mm] {}  (a2);

        \draw[-stealth, thick] (a2) --  (a3);
        \draw[-stealth, thick] (a2) to[out=-60,in=-120] node[yshift=1.5mm] {}  (a4);

        \draw[-stealth, thick] (a3) to[out=60,in=120] node[yshift=1.5mm] {} (a4);
        \draw[-stealth, thick] (a3) to[out=-60,in=-120] node[yshift=1.5mm] {}  (a4);
    \end{scope}

    \begin{scope}[xshift=210, yshift=0, scale=1]
        \vertex[fill](a1) at (1,0) {};
        \vertex[fill](a2) at (1.5,0.5) {};{};
        \vertex[fill](a3) at (2,0) {};
        \vertex[fill](a4) at (2.5,0.5) {};

        \node[] at (1,0.75) {$P$};

        \draw[thick] (a1) -- (a2);
        \draw[thick] (a3) -- (a2);
        \draw[thick] (a3) -- (a4);
    \end{scope}

    \begin{scope}[xshift=107, yshift=-60, scale=1]
        \fill[white] (0,-1) rectangle (3.5,1.5);
    
        \vertex[fill](a1) at (1,0) {};
        \vertex[fill](a2) at (1.5,0.5) {};{};
        \vertex[fill](a3) at (2,0) {};
        \vertex[fill](a4) at (2.5,0.5) {};

        \vertex[fill, orange](a5) at (1.5,-0.5) {};
        \vertex[fill, orange](a6) at (2,1) {};

        \draw[thick] (a1) -- (a2);
        \draw[thick] (a3) -- (a2);
        \draw[thick] (a3) -- (a4);

        \draw[thick, orange] (a5) -- (a1);
        \draw[thick, orange] (a5) -- (a3);
        \draw[thick, orange] (a2) -- (a6);
        \draw[thick, orange] (a4) -- (a6);

        \node[above,orange] at (a6) {\tiny $\hat{1}$};
        \node[below,orange] at (a5) {\tiny $\hat{0}$};
    \end{scope}

   \begin{scope}[xshift=107, yshift=-60, scale=1]
        \vertex[fill, cyan](a1) at (0.5,0) {};
        \vertex[fill, cyan](a2) at (1.5,0) {};
        \vertex[fill, cyan](a3) at (2,0.5) {};
        \vertex[fill, cyan](a4) at (3,0.5) {};

        
        \draw[-stealth, thick, cyan] (a1) to[out=60,in=120] node[yshift=0.75mm, xshift=5mm] {}  (a3);
        \draw[-stealth, thick, cyan] (a1) to[out=60,in=120] node[yshift=1mm,xshift=3mm] {} (a2);
        \draw[-stealth, thick, cyan] (a1) to[out=-60,in=-120] node[yshift=1.5mm] {}  (a2);

        \draw[-stealth, thick, cyan] (a2) --  (a3);
        \draw[-stealth, thick, cyan] (a2) to[out=-60,in=-120] node[yshift=1.5mm] {}  (a4);

        \draw[-stealth, thick, cyan] (a3) to[out=60,in=120] node[yshift=1.5mm] {} (a4);
        \draw[-stealth, thick, cyan] (a3) to[out=-60,in=-120] node[yshift=1.5mm] {}  (a4);
    \end{scope}

    \node[] at (7,1.8) {\small $\mathcolor{black}{G} \longmapsto \mathcolor{cyan}{P_G}$};
    \node[] at (4,-1) {\small$\mathcolor{black}{P} \longmapsto \mathcolor{cyan}{G_P}$};

\end{tikzpicture}

%% file: main.bbl
\begin{thebibliography}{10}

\bibitem{athanasiadis}
Christos~A. Athanasiadis.
\newblock Ehrhart polynomials, simplicial polytopes, magic squares and a
  conjecture of {S}tanley.
\newblock {\em J. Reine Angew. Math.}, 583:163--174, 2005.

\bibitem{baldoni-vergne}
Welleda Baldoni and Mich\`ele Vergne.
\newblock Kostant partitions functions and flow polytopes.
\newblock {\em Transform. Groups}, 13(3-4):447--469, 2008.

\bibitem{ccd}
Matthias Beck and Sinai Robins.
\newblock {\em Computing the continuous discretely}.
\newblock Undergraduate Texts in Mathematics. Springer, New York, second
  edition, 2015.
\newblock Integer-point enumeration in polyhedra, With illustrations by David
  Austin.

\bibitem{caracolvolume}
Carolina Benedetti, Rafael~S. Gonz\'{a}lez~D'Le\'{o}n, Christopher R.~H.
  Hanusa, Pamela~E. Harris, Apoorva Khare, Alejandro~H. Morales, and Martha
  Yip.
\newblock The volume of the caracol polytope.
\newblock {\em S\'{e}m. Lothar. Combin.}, 80B:Art. 87, 12, 2018.

\bibitem{combflowmodel}
Carolina Benedetti, Rafael~S. Gonz\'{a}lez~D'Le\'{o}n, Christopher R.~H.
  Hanusa, Pamela~E. Harris, Apoorva Khare, Alejandro~H. Morales, and Martha
  Yip.
\newblock A combinatorial model for computing volumes of flow polytopes.
\newblock {\em Trans. Amer. Math. Soc.}, 372(5):3369--3404, 2019.

\bibitem{berggrenpersonal}
Jonah Berggren.
\newblock Flows on gentle algebras.
\newblock Manuscript in preparation, 2024.

\bibitem{braunmcelroyfullvolumes}
Benjamin Braun and Ford McElroy.
\newblock Volume inequalities for flow polytopes of full directed acyclic
  graphs, 2024.
\newblock preprint at \url{https://arxiv.org/abs/2405.02433}.

\bibitem{BrunsHerzogCMR}
Winfried Bruns and J{\"u}rgen Herzog.
\newblock {\em Cohen-{M}acaulay rings}, volume~39 of {\em Cambridge Studies in
  Advanced Mathematics}.
\newblock Cambridge University Press, Cambridge, 1993.

\bibitem{brunsromer}
Winfried Bruns and Tim R\"{o}mer.
\newblock {$h$}-vectors of {G}orenstein polytopes.
\newblock {\em J. Combin. Theory Ser. A}, 114(1):65--76, 2007.

\bibitem{DKK}
Vladimir~I. Danilov, Alexander~V. Karzanov, and Gleb~A. Koshevoy.
\newblock Coherent fans in the space of flows in framed graphs.
\newblock In {\em 24th {I}nternational {C}onference on {F}ormal {P}ower
  {S}eries and {A}lgebraic {C}ombinatorics ({FPSAC} 2012)}, volume~AR of {\em
  Discrete Math. Theor. Comput. Sci. Proc.}, pages 481--490. Assoc. Discrete
  Math. Theor. Comput. Sci., Nancy, 2012.

\bibitem{generalpitmanstanley}
William~T. Dugan, Maura Hegarty, Alejandro~H. Morales, and Annie Raymond.
\newblock Generalized {P}itman-{S}tanley flow polytopes.
\newblock {\em S\'{e}m. Lothar. Combin.}, 89B:Art. 80, 12, 2023.

\bibitem{Ehrhart}
Eug{\`e}ne Ehrhart.
\newblock Sur les poly\`edres rationnels homoth\'etiques \`a {$n$}\ dimensions.
\newblock {\em C. R. Acad. Sci. Paris}, 254:616--618, 1962.

\bibitem{cyclicorderflow}
Rafael~S. Gonz\'{a}lez~D'Le\'{o}n, Christopher R.~H. Hanusa, Alejandro~H.
  Morales, and Martha Yip.
\newblock Column-convex matrices, {$G$}-cyclic orders, and flow polytopes.
\newblock {\em Discrete Comput. Geom.}, 70(4):1593--1631, 2023.

\bibitem{integerpointsflow}
Kabir Kapoor, Karola M\'{e}sz\'{a}ros, and Linus Setiabrata.
\newblock Counting integer points of flow polytopes.
\newblock {\em Discrete Comput. Geom.}, 66(2):723--736, 2021.

\bibitem{lidskii}
B.~V. Lidski\u{\i}.
\newblock The {K}ostant function of the system of roots {$A\sb{n}$}.
\newblock {\em Funktsional. Anal. i Prilozhen.}, 18(1):76--77, 1984.

\bibitem{gtflow}
Ricky~I. Liu, Karola M\'{e}sz\'{a}ros, and Avery St.~Dizier.
\newblock Gelfand-{T}setlin polytopes: a story of flow and order polytopes.
\newblock {\em SIAM J. Discrete Math.}, 33(4):2394--2415, 2019.

\bibitem{flowdiagonalharmonics}
Ricky~Ini Liu, Alejandro~H. Morales, and Karola M\'{e}sz\'{a}ros.
\newblock Flow polytopes and the space of diagonal harmonics.
\newblock {\em Canad. J. Math.}, 71(6):1495--1521, 2019.

\bibitem{meszaros-morales}
Karola M\'{e}sz\'{a}ros and Alejandro~H. Morales.
\newblock Volumes and {E}hrhart polynomials of flow polytopes.
\newblock {\em Math. Z.}, 293(3-4):1369--1401, 2019.

\bibitem{meszaros-morales-striker}
Karola M\'{e}sz\'{a}ros, Alejandro~H. Morales, and Jessica Striker.
\newblock On flow polytopes, order polytopes, and certain faces of the
  alternating sign matrix polytope.
\newblock {\em Discrete Comput. Geom.}, 62(1):128--163, 2019.

\bibitem{morrisidentityflow}
Alejandro~H. Morales and William Shi.
\newblock Refinements and symmetries of the {M}orris identity for volumes of
  flow polytopes.
\newblock {\em C. R. Math. Acad. Sci. Paris}, 359:823--851, 2021.

\bibitem{reinerwelker}
Victor Reiner and Volkmar Welker.
\newblock On the {C}harney-{D}avis and {N}eggers-{S}tanley conjectures.
\newblock {\em J. Combin. Theory Ser. A}, 109(2):247--280, 2005.

\bibitem{rietsch2024rootpolytopesflowpolytopes}
Konstanze Rietsch and Lauren Williams.
\newblock Root polytopes, flow polytopes, and order polytopes, 2024.
\newblock preprint at \url{https://arxiv.org/abs/2406.15803}.

\bibitem{vonbell2024triangulations}
Matias von Bell, Benjamin Braun, Kaitlin Bruegge, Derek Hanely, Zachery
  Peterson, Khrystyna Serhiyenko, and Martha Yip.
\newblock Triangulations of flow polytopes, ample framings, and gentle
  algebras.
\newblock {\em Selecta Math. (N.S.)}, 30(3):Paper No. 55, 2024.

\bibitem{framedtriangulations}
Matias von Bell, Rafael~S. Gonz\'{a}lez~D'Le\'{o}n, Francisco~A.
  Mayorga~Cetina, and Martha Yip.
\newblock On framed triangulations of flow polytopes, the {$\nu$}-{T}amari
  lattice and {Y}oung's lattice.
\newblock {\em S\'{e}m. Lothar. Combin.}, 85B:Art. 42, 12, 2021.

\end{thebibliography}
